\documentclass[12pt,english]{article}
\usepackage[cp1251]{inputenc}
\usepackage{geometry}
\usepackage{color}
\usepackage{amsmath}
\usepackage{amssymb}
\usepackage{esint}
\usepackage{amsthm}

\def\E{\hskip.15ex\mathrm{E}\hskip.10ex}
\def\P{\mathrm{P}}

\def\phi{\varphi}



\usepackage{babel}

\newtheorem{theorem}{Theorem}
\newtheorem{lemma}[theorem]{Lemma}
\newtheorem{remark}[theorem]{Remark}
\newtheorem{definition}[theorem]{Definition}
\newtheorem{exercise}[theorem]{Exercise}
\newtheorem{Assumption}[theorem]{Assumption}
\newtheorem{corollary}[theorem]{Corollary}
\newtheorem{proposition}[theorem]{Proposition}
\newtheorem{openq}[theorem]{Open Question}

\begin{document}
\global\long\def\E{\mathbb{E}}
\global\long\def\P{\mathbb{P}}
\global\long\def\N{\mathbb{N}}
\global\long\def\ind{\mathbb{I}}

\title{
{\normalsize\tt\hfill\jobname.tex}\\
Ergodic Markov processes and Poisson equations\\ (lecture notes)}
\author{A. Yu. Veretennikov\footnote{ University of Leeds, UK; National Research University Higher School of Economics, and Institute for Information Transmission Problems, Moscow, Russia, email: a.veretennikov @ leeds.ac.uk. The work was prepared within the framework of a subsidy granted to the HSE by the Government of the Russian Federation for the implementation of the Global Competitiveness Program, and supported by the RFBR grant 14-01-00319-a. The author is also grateful to two anonymous referees for many valuable remarks.}
}

\maketitle
\begin{abstract}
These are lecture notes on the subject defined in the title. As such, they do not pretend to be really new, perhaps, except for the section  \ref{sec:pp} about Poisson equations with potentials; also, the 
convergence rate shown in (\ref{newrate}) -- (\ref{newrate2}) is possibly less-known. Yet, the hope of the author is that these notes may serve as a bridge to the important area of Poisson equations ``in the whole space'' and with a parameter, the latter theme not being presented here. Why this area is so important was explained in many papers and books including \cite{EthierKurtz, PSV,  PV03}: it provides one of the main tools in diffusion approximation in the area stochastic averaging. Hence, the aim of these lectures is to prepare the reader to ``real'' Poisson equations -- i.e., for differential operators instead of difference operators -- and, indeed, to diffusion approximation. Among other presented topics  we mention coupling method and convergence rates in the Ergodic theorem.

\end{abstract}

\tableofcontents

\section{Introduction} 
In these lecture notes we will consider the following issues: Ergodic theorem (in some textbooks called Convergence theorem, while Ergodic would be reserved for what we call Law of Large Numbers -- see below), Law of Large Numbers (LLN), Central Limit Theorem (CLT), Large Deviations (LDs) for Markov chains (MC), and as one  of the most important applications, a Poisson equation. LLN, CLT and LDs are the basis of most of  statistical applications. Everything is presented on the simplest model of a Markov chain with positive transition probabilities on a finite state space, and in some cases we show a few more general results where it does not require too much of additional efforts. This simplified version may be regarded as a preparation to more advanced situations of Markov chains on more general state space, including non-compact ones and including Markov diffusions. A special place in this plan is occupied by coupling method, a famous idea, which is not necessary for any result in these lecture notes; yet, it is a rather convenient tool ``for thinking'', although sometimes not very easy for a rigorous presentation. We show the Ergodic theorem firstly without and then with coupling method. Poisson equations in this paper are discrete analogues of ``real'' Poisson equations for elliptic differential operators of the second order in mathematical physics. We consider equations {\em without} a potential -- the most useful tool in diffusion approximations, cf. \cite{EthierKurtz,  PSV, PV03} -- and also {\em with} a potential. The problem of smoothness of solutions with respect to a parameter -- which makes this stuff so important in diffusion approximations and which is one of the main motivations of the whole theory -- is not presented; however, these notes may be regarded as a bridge to this smoothness issue.

~

These notes are based on several different courses delivered by the author at various universities in various years, including Moscow State University, Helsinki Technical University (now Aalto university),  University of Leeds and Magic consortium (http://maths-magic.ac.uk/index.php), and National Research University Higher School of Economics -- Moscow. The author thanks all participants -- not only students -- for their interest and patience and for many useful remarks.

The initial plan involved non-compact cases with problems related to stability or recurrence properties of processes in such spaces. However, this would require significantly more time and many more pages. Hence, this more ambitious task is postponed for some future.

Some classical results are given without proofs although they were proved in the courses delivered. The references on all such ``missing'' proofs are straightforward.

Finally, let us mention that the following numeration system is accepted here: all items such as Theorem, Lemma, Remark, Definition and some others are numbered by a unique sequence of natural numbers. This method was accepted in some well-known textbooks and the author shares the view about its convenience.

~

\noindent
The following notations will be used for a process $(X_n, \, n\ge 0)$: 
$$
{\cal F}^X_n = \sigma(X_k: \, k\le n); \quad {\cal F}^X_{(n)} = \sigma(X_n).
$$ 
The following notations from the theory of Markov processes will be accepted (cf. \cite{Dynkin}): the index $x$ in $\mathbb E_x$ or $\mathbb P_x$ signifies the expectation or the probability measure related to the non-random initial state of the process $X_0$. This initial state may be also random with some distribution $\mu$, in which case notations  $\mathbb E_\mu$ or $\mathbb P_\mu$ may be used. 

If state space $S$ is finite, then $|S|$ denotes the number of its elements. In the sequel $\cal P$ denotes the transition matrix $\left(p_{ij}\right)_{1\le i,j \le |S|}$ of the process in the cases where state space of the process is finite.

Since this is a course about ergodic properties, we do not recall the definitions of what are Markov, strong Markov, homogeneous Markov processes (MP) which are assumed to be known to the reader: consult any of the textbooks \cite{Bor98, Doob53, Dynkin, EthierKurtz, Karlin,  Kry-ln, Seneta1, Wentzell} if in doubt.

\section{Ergodic Theorem -- 1}\label{sec1}
In this section we state and prove a simple ergodic theorem for Markov chains on a finite state space. However, we start with a more general setting because later in the end of these lecture  notes a more general setting will be addressed. 
Ergodic Theorem for Markov chains in a simple situation of finite state spaces is due to Markov,  although sometimes it is attributed to Kolmogorov with a reference to Gnedenko's textbook, and sometimes to Doeblin (see \cite{Gnedenko}, \cite{Doeblin}). We emphasize that this approach was introduced by Markov himself (see \cite{Markov, Seneta1, Seneta2}). Kolmogorov, indeed, has contributed to this area: see, in particular, \cite{Kolm49}.

Let us consider a homogeneous Markov chain $X = (X_n), \, n=0,1,2,\ldots$ with a general topological state space $(S, {\cal S})$ where ${\cal S}$ is the family of all Borel sets in $S$ assuming that ${\cal S}$ contains all single point subsets. Let $P_x(A)$ be the transition kernel, that is, $P_x(A) = \mathbb P(X_1\in A | X_0=x)  \equiv \mathbb P_x(X_1\in A) $; recall that for any $A \in {\cal S}$ this function is assumed Borel measurable in $x$ (see \cite{Dynkin}) and a measure in $A$ (of course, for a finite $S$ this is not a restriction). Denote by $P_x(n,A)$ the $n$-step transition kernel, i.e., $P_x(n,A) = \mathbb P_x(X_n\in A)$; for a finite Markov chain and if $A=j$, the notation $p_{ij}^{(n)}$ will be used, too. If initial state is random with distribution $\mu$, we will be using a similar notation $P_\mu(n,A)$ for the probability $\mathbb P_\mu(X_n\in A)$. Repeat that $\mathbb P_{inv}(X_n\in A)$ signifies $P_\mu(X_n\in A)$ with the (unique) invariant measure $\mu$; naturally, this value does not depend on~$n$. 

~

Recall the definition of ergodicity for Markov chains (MC). 

\begin{definition}
An MC  $(X_n)$ is called Markov ergodic iff the sequence of transition measures $(\mathbb P_x(n, \cdot))$ has a limit in total variation metric, which is a probability measure and if, in addition, this limiting measure does not depend on $x$, 
\begin{equation}\label{erg-def}
 \lim_{n\to\infty} \mathbb P_x(n,A) = \mu(A), \quad \forall \;  A\in {\cal S}.
\end{equation}
\end{definition}
\noindent
Recall that the total variation distance or metric between two probability measures may be defined as 
\[
\|\mu-\nu\|_{TV} := 2 \sup\limits_{A \in {\cal S}} (\mu(A) - \nu(A)). 
\]

\begin{definition}\label{irreduci}
An MC $(X_n)$ is called 
irreducible iff for any $x\in S$ and $A\in {\cal S}, \, A \not = \emptyset$, there exists $n$ such that 
\[
\mathbb P_x(X_n \in A) > 0. 
\]
An MC $(X_n)$ is called $\nu$-irreducible for a given measure $\nu$ on $(S, {\cal S})$ iff for any $x\in S$ and $A\in {\cal S}, \, \nu(A)>0$ there exists $n$ such that 
\[
\mathbb P_x(X_n \in A) > 0. 
\]
\end{definition}
Of course, weaker or stronger ergodicity properties (definitions) may be stated with weaker, or, respectively, stronger metrics. Yet, in the finite state space case most of them -- although, not all of them -- are  equivalent.

\begin{exercise}
In the case of a finite state space $S$ with ${\cal S}= 2^S$ (all subsets of $S$) and a {\em counting} measure \(\nu\) such that $\nu(A) = |A| := \mbox{the number of elements in $A \subset S$}$, show that 
$\nu$-irreducibility of a MC is equivalent to the claim that there exists $n>0$ such that the $n$-step transition probability matrix ${\cal P}^n$ is {\bf positive}, that is, all elements of it are strictly positive.  
\end{exercise}
\noindent
The most standard is the notion of $\nu$-irreducibility of an MC where $\nu$ is the unique invariant measure of the process. 

\begin{definition}
Stationary or invariant probability measure $\mu$ for a Markov process $X$ is a measure on $\cal S$ such that for each $A\in \cal S$ and any $n$, 
$$
\mu(A) = \sum_{x\in  S} \mu(x) P_x(n,A). 
$$
\end{definition}

\begin{lemma}\label{sta-condition}
A probability measure $\mu$ is stationary for $X$ iff 
$$
\mu \cal P = \mu, 
$$
where $\cal P$ is the transition probability matrix of $X$. 
\end{lemma}
\noindent
{\em Proof} is straightforward by induction. 

\begin{lemma}\label{finite-sta}
For any (homogeneous) Markov chain in a finite state space $S$ there exists at least one stationary measure. 
\end{lemma}
\noindent
{\it Proof of the Lemma \ref{finite-sta}.} The method is due to Krylov and Bogoliubov (Kryloff and Bogoliuboff,  \cite{Kryloff}). Let us fix some (any) $i_0\in S$, and consider Ces\`aro averages
\[
\frac1{n+1}\sum_{k=0}^{n}p_{i_0,j}^{(k)}, \quad 1\le j \le N, \quad n\ge 1,
\]
where $N=|S|$. 
Due to the boundedness, this sequence of vectors as $n\to\infty$ has a limit over some subsequence, say, $n'\to\infty$,
\[
\frac1{n'+1}\sum_{k=0}^{n'}p_{i_0,j}^{(k)} \to \pi_j, \quad 1\le j \le N, \quad n'\to 1, 
\]
where by the standard convention, $p^{(0)}_{ij} = \delta_{ij}$ (Kronecker's symbol). 
Since $S$ is finite, it follows that $(\pi_j, \, 1\le j\le N)$ is a probability distribution on $S$. Finally, stationarity follows from the following calculus based on Chapman--Kolmogorov's equations, 
\begin{eqnarray*}
&\displaystyle \frac1{n'+1}\sum_{k=0}^{n'}p_{i_0,j}^{(k)} = 
\frac1{n'+1}\sum_{k=0}^{n'} \sum_{\ell = 1}^{N} 
p_{i_0,\ell}^{(k-1)} p_{\ell, j}^{} +  \frac1{n'+1} p_{i_0,j}^{(0)}
 \\\\
&\displaystyle =  \sum_{\ell = 1}^{N} \frac1{n'+1}\sum_{k=0}^{n'-1} 
p_{i_0,\ell}^{(k)} p_{\ell, j}^{} +  \frac1{n'+1} p_{i_0,j}^{(0)} 
 \\\\
&\displaystyle =   \sum_{\ell = 1}^{N} \frac1{n'+1}\sum_{k=0}^{n'} 
p_{i_0,\ell}^{(k)} p_{\ell, j}^{} +  \frac1{n'+1} p_{i_0,j}^{(0)} 
- \frac1{n'+1} \sum_{\ell = 1}^{N} 
p_{i_0,\ell}^{(n')} p_{\ell, j}^{}.
\end{eqnarray*}
It follows, 
\begin{eqnarray*}
\lim_{n' \to \infty} \frac1{n'+1}\sum_{k=0}^{n'}p_{i_0,j}^{(k)} 
= \sum_{\ell = 1}^{N} 
\pi_\ell p_{\ell, j}^{}.
\end{eqnarray*}
Hence, 
\begin{eqnarray*}
\pi_j = \sum_{\ell = 1}^{N} \pi_j p_{\ell, j}^{} \quad \sim \quad \pi = \pi \cal P.  
\end{eqnarray*}
Hence, the distribution $(\pi_j)$ is stationary due to the Lemma \ref{sta-condition}. The Lemma \ref{finite-sta} is proved. 

~
\begin{remark}\label{rem:brouwer}
Note that for a finite $S$ the statement of the Lemma, actually, may be proved much faster by applying the Brouwer fixed point theorem, as it is done, for example, in \cite{tutu}. Yet, the method used in the proof seems deeper, and it can be used in a much more general situation including ``non-compact'' cases. (However, we are not saying that the use of Brouwer's fixed point theorem is restricted to finite state spaces.)
\end{remark}

~

{\bf From now on, in this and several following sections we consider the case of a finite state space $S$; a more general case will be addressed in the last sections.} 
The next condition suggested by Markov himself plays a very important role in the analysis of asymptotic behaviour of a (homogeneous) Markov chain (MC in the sequel).  Let there exist $n_0$ such that 
\begin{equation}\label{Mar}
 \kappa_{n_0} := \inf_{i,i'}  \sum_j \min(P_i(n_0,j), P_{i'}(n_0,j)) 
 \equiv \inf_{i,i'}  \sum_j \min(p^{n_0}_{i,j}, p^{n_0}_{i',j}) > 0. 
\end{equation}
By the suggestion of S. Seneta, this coefficient  $\kappa_{n_0}$ (as well as $\kappa_{}$ in (\ref{Mar2}) and 
in (\ref{MD})) is properly called Markov--Dobrushin's.

\noindent
Unlike in the continuous time case, in discrete time situation there are potential complications related to possible {\em cycles}, that is, to a {\em periodic structure} of the process. A typical example of such a periodic structure is a situation where the state space is split into two parts, $S =  S_1 \cup S_2$, which do not intersect, and $X_{2n} \in S_1$, while $X_{2n+1} \in  S_2$ for each $n$. Then ergodic properties is reasonable to study separately for $Y_n:= X_{2n}$ and  for $Z_n:= X_{2n+1}$. In other words, this complication due to periodicity does not introduce any real news, and by this reason there is a tradition to avoid this situation. Hence, in the sequel we will study our ergodic process under the assumption $n_0=1$ in the condition (\ref{Mar}). Similar results could be obtained under a more general assumption of {\em aperiodicity}.

So, here is the simplified version of (\ref{Mar}), which will be accepted in the sequel: 
\begin{equation}\label{Mar2}
 \kappa := \inf_{i,i'}  \sum_j \min(P_i(1,j), P_{i'}(1,j)) 
 \equiv \inf_{i,i'}  \sum_j \min(p^{}_{ij}, p^{}_{i'j}) > 0. 
\end{equation}
\noindent
Also, to clarify the ideas we will be using  in some cases the following stronger assumption, 
\begin{equation}\label{Mar3}
\kappa_0:=\inf_{ij} p_{ij} > 0.
\end{equation}
However, eventually, the assumption (\ref{Mar3}) will be dropped and only (\ref{Mar2}) will remain in use.

\begin{theorem}\label{thm_erg1}
Let the assumption (\ref{Mar2}) hold true. Then the process $(X_n)$ is ergodic, i.e., there exists  a limiting probability measure $\mu$ such that (\ref{erg-def}) holds true. Moreover, the uniform bound is satisfied for every $n$, 
\begin{equation}\label{exp_bd}
 \sup_{x}\sup_{A\in {\cal S}} |P_x(n,A) - \mu(A)| \le (1-\kappa)^{n},   
\end{equation}
and the measure $\mu$ is a unique invariant one.
\end{theorem}
\noindent
{\em Proof of Theorem \ref{thm_erg1}} is classical and may be found in many places, for example, in \cite{Gnedenko}. 

~

\noindent
(A) Denote for any $A$,
$$
m^{(n)}(A):= \min_i P_i(n,A), \quad M^{(n)}(A):= \max_i P_i(n,A).
$$
By Chapman--Kolmogorov's equation,
\begin{eqnarray*}
& \displaystyle m^{(n+1)}(A)= \min_i P_i(n+1,A) = \min_i \sum_j p_{ij} P_j(n,A)
 \\ \\
& \displaystyle \ge  \min_i \sum_j p_{ij} \min_{j'} P_{j'}(n,A) =m^{(n)}(A),
\end{eqnarray*}
which signifies that the sequence $m^{(n)}(A)$ does not decrease in $n$. Similarly, the sequence $M^{(n)}(A)$ does not increase in $n$. Hence, it suffices to show that
\begin{equation}\label{meq}
M^{(n)}(A) - m^{(n)}(A) \le (1-\kappa)^n.
\end{equation}
(B) Again by Chapman--Kolmogorov's equation,
\begin{eqnarray*}
M^{(n)}(A) - m^{(n)}(A)= \max_i P_i(n,A) -
\min_{i'} P_{i'}(n,A)
 \\ \\
=  \max_i \sum_j p_{ij} P_{j}(n-1,A) -
\min_{i'} \sum_j p_{i' j} P_{j}(n-1,A).
\end{eqnarray*}
Let maximum here be attained at $i_+$ while minimum at $i_-$. Then,
\begin{eqnarray}\label{30}
& \displaystyle M^{(n)}(A) - m^{(n)}(A)= \sum_j p_{i_+ j} P_{j}(n-1,A) - \sum_j p_{i_- j} P_{j}(n-1,A)
 \nonumber \\ \nonumber  \\
& \displaystyle =  \sum_j (p_{i_+ j} - p_{i_- j}) P_{j}(n-1,A).
\end{eqnarray}
(C) Denote by $S^+$ the part of the sum in the right hand side of (\ref{30}) with just $(p_{i_+ j} - p_{i_- j})\ge 0$, and by $S^-$ the part of the sum with $(p_{i_+ j} - p_{i_- j}) < 0$. Using notations $a_+ = a \vee 0$ and $a_- = a\wedge 0$ (where $a\vee b = \max(a,b)$ and $a\wedge b = \min(a,b)$), we estimate,
$$
S^+ \le \sum_j (p_{i_+ j} - p_{i_- j})_+ M^{(n-1)}(A) 
= M^{(n-1)}(A)\sum_j (p_{i_+ j} - p_{i_- j})_+,
$$
and
$$
S^- \le  \sum_j (p_{i_+ j} - p_{i_- j})_- m^{(n-1)}(A).
$$
Therefore,
\begin{eqnarray*}
& \displaystyle M^{(n)}(A) - m^{(n)}(A)= S^+ + S^-
 \\ \\
& \displaystyle \le M^{(n-1)}(A) \sum_j (p_{i_+ j} - p_{i_- j})_+ +  m^{(n-1)}(A) \sum_j (p_{i_+ j} - p_{i_- j})_-.
\end{eqnarray*}
(D) It remains to notice that
\[
\sum_j (p_{i_+ j} - p_{i_- j})_- = - \sum_j (p_{i_+ j} - p_{i_- j})_+,
\]
and
\begin{equation}\label{32}
\sum_j (p_{i_+ j} - p_{i_- j})_+ \le 1-\kappa.
\end{equation}
The first follows from the normalization condition
$$
\sum_j p_{i_+ j} = \sum_j p_{i_- j} = 1,
$$
while the second from (recall that $(a-b)_+ = a - a\wedge b \equiv \min(a,b)$ for any real values $a, b$)
\begin{eqnarray*}
& \displaystyle \sum_j (p_{i_+ j} - p_{i_- j})_+
= \sum_j (p_{i_+ j} - \min(p_{i_- j},p_{i_+ j}))
 \\\\
& \displaystyle = 1  -  \sum_j \min(p_{i_- j},p_{i_+ j})
\le 1-\kappa
\end{eqnarray*}
(see the definition of $\kappa$ in (\ref{Mar2})). So, we find that
$$
M^{(n)}(A) - m^{(n)}(A) \le (1-\kappa)\,(M^{(n-1)}(A) - m^{(n-1)}(A)).
$$
By induction this implies (\ref{meq}). So, (\ref{exp_bd}) and  uniqueness of the limits $\pi_j=\lim_{n\to\infty} p^{(n)}_{ij}$ follow. 

~

\noindent
(E) The invariance of the measure $\mu$ and uniqueness of the invariant measure follow, in turn, from (\ref{exp_bd}). Indeed, let us start the process from any invariant distribution $\mu$ -- which exists due to the Lemma \ref{finite-sta} --  then $\mu_j\equiv \mathbb P_\mu(X_n=j)= \sum_{\ell}\mu_\ell p^{(n)}_{ij} \to \pi_j, \; n\to\infty$. However, the left hand side here does not depend on $n$. Hence, $\mu_j=\pi_j$. The Theorem \ref{thm_erg1} is proved.

~

\noindent
Recall that the {\bf total variation distance} or metric between two probability measures may be defined as 
\[
\|\mu-\nu\|_{TV} := 2 \sup_A (\mu(A) - \nu(A)). 
\]
Hence, the inequality (\ref{exp_bd}) may be rewritten 
as 
\begin{equation}\label{exp_bd2}
 \sup_{x} \|P_x(n,\cdot) - \mu(\cdot)\|_{TV} \le 2 (1-\kappa)^{n}. 
\end{equation}

\begin{corollary}\label{cor6}
Under the assumption of the Theorem \ref{thm_erg1}, for any bounded Borel function $f$ and for any $0\le s < t$, 
\[
\sup_x |{\mathbb E}_x (f(X_t) | X_s) - {\mathbb E}_{inv} f(X_t)| 
\equiv \sup_x |{\mathbb E}_x (f(X_t) - {\mathbb E}_{inv} f(X_t)| X_s) | 
\le C_f (1-\kappa)^{t-s}, 
\]
or, equivalently, 
\[
\sup_x |{\mathbb E}_x (f(X_t) | {\cal F}^X_s) - {\mathbb E}_{inv} f(X_t)| 
\le C_f (1-\kappa)^{t-s}, 
\]
where $C_f = \max\limits_j |f(j)| \equiv \| f\|_{B(S)}$.
\end{corollary}
\noindent
This useful Corollary follows from the Theorem \ref{thm_erg1}.

~

It is worth noting that in a general case there is a significantly weaker condition than (\ref{Mar}) (or, in the general case weaker than (\ref{MD}) -- see below in the section \ref{sec:erg_gen}), which also guarantees an exponential convergence rate to a unique invariant measure. We will show this condition -- called Doeblin-Doob's one -- and state the corresponding famous Doeblin--Doob's theorem on convergence, but for the proof we refer the reader to \cite{Doob53}. 

\begin{definition}[DD-condition] There exist a finite (sigma-additive) measure $\nu\ge 0$ and $\epsilon>0$, $s>0$ such that $\nu(A)\le \epsilon$ implies 
\[
\sup_x P_x(s, A) \le 1 - \epsilon. 
\]
\end{definition}

\begin{theorem}[Doeblin--Doob, without proof]
If the DD-condition is satisfied for an {\em aperiodic MP with a unique class of ergodicity (see \cite{Doob53})} on the state space $\cal S$, then there exist $C,c>0$ such that 
\begin{equation}\label{exp_bd2}
 \sup_{x}\sup_{A\in {\cal S}} |P_x(n,A) - \mu(A)| \le C\exp(-cn), \quad n\ge 0. 
\end{equation}
\end{theorem}
It turns out that under the assumption (DD), the constants in the upper bound (\ref{exp_bd2}) {\em cannot be effectively computed}, i.e., they may be arbitrary even for the same $\epsilon$ and $\nu$, say. This situation dramatically differs from the case of conditions (\ref{Mar3}) and (\ref{Mar2}), where both  constants in the upper bound are effectively and explicitly evaluated.

\begin{openq}
It is interesting whether or not there may exist any intermediate situation with a bound like (\ref{exp_bd2}) -- in particular, it should be uniform in the initial state -- with {\em computable} constants $C,c$ under an assumption lying somewhere in ``between'' Markov--Dobrushin's and Doeblin--Doob's. 
Apparently, such a condition may be artificially constructed from a ``non--compact'' theory with an exponential recurrence, but then the bounds would not be uniform in the initial data. In fact, some relatively simple version of a desired condition will be shown in the end of this text, see the Theorem \ref{lastthm}. However, it does not totally close the problem, e.g., for non-compact spaces. 
\end{openq}

\section{LLN for homogeneous MC, finite $S$}
It may seem as if the Ergodic Theorem with uniform exponential convergence rate in total variation metric were all we could wish about ergodic features of the process. Yet, the statement of this theorem itself even does not include the Law of Large Numbers (LLN), which is not emphasized in most of the textbooks. However, the situation with LLN (as well as with Central Limit Theorem -- CLT) is good enough, which is demonstrated below. The Theorem \ref{thm_lln1} under the assumption (\ref{Mar3}) belongs to A.A. Markov, see \cite{Markov, Seneta1}. 

\begin{theorem}[Weak LLN]\label{thm_lln1}
Under the assumptions of the Theorem \ref{thm_erg1}, for any function $f$ on a finite state space $S$, 
\begin{equation}\label{eq_lln1}
\frac1n \sum_{k=0}^{n-1} f(X_k) \stackrel{\mathbb P}{\to} {\mathbb E}_{inv}f(X_0), 
\end{equation}
where ${\mathbb E}_{inv}$ stands for the expectation of $f(X_0)$ with respect to the invariant probability measure of the process, while $\mathbb P$ denotes the measure, which corresponds to the initial value or distribution of $X_0$: the latter may be, or may be not stationary. 
\end{theorem}
\noindent
NB. Note that a simultaneous use of stationary and non-stationary measures is not a contradiction here. The initial state could be either non-random, or it may have some distribution. At the same time, the process has a unique invariant measure, and the writing $\mathbb E_{inv}f(X_0)=0$ signifies the mere fact that $\sum\limits_{y\in S} f(y) \mu(y) = 0$, but it is in no way in a conflict with a non-stationary  initial distribution. In the next proof we use $\mathbb P$ and,  respectively, $\mathbb E$ without specifying the initial state or distribution. However, this initial distribution (possibly concentrated at one single state) exists and it is fixed throughout the proof. 

~

\noindent
{\em Proof of the Theorem \ref{thm_lln1}.} {\bf 1}. First of all, note that (\ref{eq_lln1}) is equivalent to
\[
\frac1n \sum_{k=0}^{n-1} (f(X_k) - {\mathbb E}_{inv}f(X_0)) \stackrel{\mathbb P}{\to} 0, 
\]
so, without loss of generality we may and will assume that $ {\mathbb E}_{inv}f(X_0) = 0$. Now we estimate with any $\epsilon > 0$ by the Bienaym\'e--Chebyshev--Markov inequality, 
\begin{eqnarray}
{\mathbb P} \left(|\frac1n \sum_{k=0}^{n-1} f(X_k)| > \epsilon\right)
\le \epsilon^{-2}n^{-2}{\mathbb E} |\sum_{k=0}^{n-1} f(X_k)|^2
 \nonumber \\ \label{ett}\\ \nonumber
=  \epsilon^{-2}n^{-2}{\mathbb E} \sum_{k=0}^{n-1} f^2(X_k)
+ 2\epsilon^{-2}n^{-2}{\mathbb E} \sum_{0\le k<j\le n-1}^{} f(X_k)f(X_j).
\end{eqnarray}
Here the first term, clearly (as $f$ is bounded), satisfies, 
\[
  \epsilon^{-2}n^{-2}{\mathbb E} \sum_{k=0}^{n-1} f^2(X_k) \to 0, \quad n\to\infty.
\]
Let us transform the second term as follows for $k<j$: 
\begin{eqnarray*}
{\mathbb E} f(X_k)f(X_j) = 
{\mathbb E} (f(X_k) {\mathbb E} (f(X_j) | X_k)),  
\end{eqnarray*}
and recall that due to the Corollary \ref{cor6} to the Ergodic theorem, 
\[
|{\mathbb E} (f(X_j) | X_k) - {\mathbb E}_{inv} f(X_j)| \le C_f (1-\kappa)^{j-k}, 
\]
where due to our convention ${\mathbb E}_{inv} f(X_j)=0$.
Therefore, we have, 
\begin{eqnarray*}
&\displaystyle |{\mathbb E} \sum_{k<j}^{} f(X_k)f(X_j)| = |{\mathbb E} \sum_{k<j}^{} f(X_k) {\mathbb E} (f(X_j) | X_k)| 
 \\\\
&\displaystyle \le C_f \sum_{k,j:\, 0\le k<j<n}^{}(1-\kappa)^{j-k} \le C n, \quad \mbox{with} \;\; C = C_f \kappa^{-1}. 
\end{eqnarray*}
Thus, the second term in (\ref{ett}) also goes to zero as $n\to\infty$. The Theorem \ref{thm_lln1} is proved. 

\begin{remark}
Recall that $f$ is bounded and exponential rate of convergence is guaranteed by the assumptions. This suffices for a strong LLN via higher moments for sums. However, it will not be used in the sequel, so we do not show it here. 
\end{remark}

\section{CLT, finite $S$}
In this section state space $S$ is also finite. For the function $f$ on $S$, let 
\begin{equation}\label{sig}
\sigma^2:= {\mathbb E}_{inv}(f(X_0) -\mathbb E_{inv}f(X_0))^2 
+ 2 \sum_{k=1}^{\infty} {\mathbb E}_{inv}(f(X_0) -\mathbb E_{inv}f(X_0))(f(X_k) - \mathbb E_{inv}f(X_k)). 
\end{equation}
It is known that this definition provides a non-negative value (for completeness, see the two lemmata below). 

\begin{lemma}\label{sge0}
Under our standing assumptions ($S$ is finite  and $\min\limits_{ij}p_{ij}>0$), 
\begin{equation}\label{sp}
\sigma^2\ge 0, 
\end{equation}
and, moreover, 
\begin{eqnarray}\label{varn0}
n^{-1}{\mathbb E}_{inv} \left(\sum^{n-1}_{r=0} \left(f(X_r) - \mathbb E_{inv}f(X_0)\right)\right)^2 \to  \sigma^2, \quad n\to\infty, 
\end{eqnarray}
\noindent
where the latter convergence is irrespectively of whether $\sigma^2>0$, or $\sigma^2=0$.
\end{lemma}

{\em Proof.} 
Without loss of generality, we may and will assume now that  ${\mathbb E}_{inv}f(X_0) =0$ (otherwise, this mean value can be subtracted from $f$ as in the formula (\ref{varn0})). 
Note also that in this case the variance of the random variable $\displaystyle n^{-1/2}\sum\limits^{n-1}_{r=0} f(X_r)$ computed {\bf with respect to the invariant measure} coincides in this case with its second moment. Since ${\mathbb E}_{inv}f(X_i) = 0$ for any $i$,  this second moment may be evaluated as follows, 
\begin{eqnarray*}
&\displaystyle n^{-1} {\mathbb E}_{inv}(\sum^{n-1}_{r=0} f(X_r))^2
= {\mathbb E}_{inv}f^2(X_0) 
+ 2n^{-1} \sum_{0\le i < j\le n-1} {\mathbb E}_{inv}f(X_i)f(X_j) 
 \\\\
&\displaystyle =  {\mathbb E}_{inv}f^2(X_0) 
+ 2n^{-1} \sum_{r=1}^{n-1} (n-r) {\mathbb E}_{inv}f(X_0)f(X_r) 
 \\\\
&\displaystyle \stackrel{\mbox{\tiny clearly}}{\to} {\mathbb E}_{inv}f^2(X_0) 
+ 2 \sum_{r=1}^{\infty} {\mathbb E}_{inv}f(X_0)f(X_k)= \sigma^2, \quad n\to\infty. 
\end{eqnarray*}
Here the left hand side is non-negative, so $\sigma^2$ is non-negative, too. The Lemma \ref{sge0} is proved.

\begin{lemma}\label{s<infty}
Under the same assumptions as in the previous Lemma, $\sigma^2 < \infty.$
\end{lemma}
\noindent
{\em Proof.} Again, without loss of generality, we may and will assume $\bar f:= {\mathbb E}_{inv}f(X_0) =0$,  
and $\|f\|_{B}\le 1$. We have, due to the Corollary \ref{cor6} applied with $\bar f = 0$,
\begin{eqnarray*}
|{\mathbb E}_{inv} f(X_0)f(X_k) | = 
|{\mathbb E}_{inv} f(X_0){\mathbb E}_{inv}(f(X_k) | X_0)|
\le  C_f{\mathbb E}_{inv} |f(X_0)| q^{k}, 
\end{eqnarray*}
with some $0\le q < 1$ and $C_f = \|f\|_{B} \le 1$. 
So, the series in (\ref{sig}) does converge and the Lemma \ref{s<infty} is proved.

\begin{theorem}\label{thm_clt1}
Let the assumption (\ref{Mar2}) hold true. Then for any function  $f$ on $S$, 
\begin{equation}\label{eq_clt1}
\frac1{\sqrt{n}} \sum_{k=0}^{n-1} (f(X_k) - {\mathbb E}_{inv}f(X_0)) \stackrel{\mathbb P}{\Longrightarrow} \eta \sim {\cal N}(0, \sigma^2),
\end{equation} 
where $\Longrightarrow$ stands for the weak convergence with respect to the original probability measure (i.e., generally speaking, non-invariant). 
\end{theorem}
Emphasize that we subtract the expectation with respect to the invariant measure, while weak convergence holds true with respect to the initial measure, which is not necessarily invariant. (We could have subtract the actual expectation instead; the difference would have been negligible due to the Corollary \ref{cor6}.)

\begin{remark}
About Markov's method in CLT the reader  may consult the textbook \cite{tutu}. Various approaches  can be found in 
\cite{bernstein26, bernstein44, Doob53,  Kolm49, Nagaev68, Nagaev15, Seneta1}, et al. For a historical review see \cite{Seneta2}. 
A nontrivial issue of distinguishing the cases  $\sigma^2>0$ and  $\sigma^2=0$ for stationary Markov chains is under discussion in \cite{Bezhaeva} for finite MC where a criterion has been established for $\sigma^2=0$; this criterion was extended to more general cases in \cite{KM}. A simple example of irreducible aperiodic MC (with $\min\limits_{ij}p_{ij}=0$) and a non-constant function $f$ where  $\sigma^2=0$ can be found in \cite[ch. 6]{tutu}. Nevertheless, there is a general belief that ``normally'' in ``most of cases'' $\sigma^2>0$. 
(Recall that zero (a constant) is regarded as a degenerate Gaussian random variable ${\cal N}(0,0)$.) On using weaker norms in CLT for Markov chains see~\cite{Kulik}. 
\end{remark}

~

\noindent
{\em Proof of the Theorem \ref{thm_clt1}.} 
Without loss of generality, assume that $\|f\|_B\le 1$, and that 
$$
{\mathbb E}_{inv}f(X_0) = 0. 
$$

~

\noindent
{\bf I.} Firstly, consider the case $\sigma^2 > 0$. We want to check the assertion, 
\[
\mathbb E\exp(i\frac{\lambda}{\sqrt{n}}\sum_{r=0}^{n} f(X_r)) 
\to \exp(-\lambda^2 \sigma^2/2), \quad n\to\infty.
\]
In the calculus below there will be expactations with respect to the measure $\mathbb P$ (denoted by $\mathbb E$) and some other expectations $\mathbb E_{inv}$. Note that they are different: the second one means expectation of a function of a random variable $X_k$ computed with respect to the invariant measure of this process. 

\noindent
We are going to use Bernstein's method of ``corridors and windows'' (cf. \cite{bernstein26, bernstein44}). 
Let us split the interval $[0,n]$ into partitions of two types: larger ones called ``corridors'' and smaller ones called ``windows''. Their sizes will increase with $n$ as follows. Let $\displaystyle k: = [n/[n^{3/4}]]$ be the total number of long corridors of equal length (here $[a]$ is the integer part of $a\in \mathbb R$); this length will be chosen shortly as equivalent to $n^{3/4}$. The length of each window is $w:= [n^{1/5}]$. Now, the length of each corridor except the last one is $c:= [n/k] - w \equiv [n/k] - [n^{1/5}]$; the last complementary corridor has the length $c_{L} := n - k[n/k]$; note that $c_{L} \le [n/k] \sim n^{3/4}$, $k \sim n^{1/4}$ (i.e., $k/n^{1/4} \to 1, \, n\to \infty$), and $c \sim n^{3/4}$. 

The total length of all windows is then equivalent to $w \times k \sim n^{1/5 + 1/4} = n^{9/20}$; note for the sequel that $n^{9/20} <\!\!\!< n^{1/2}$. As was mentioned earlier, the length of the last corridor does not exceed $k$, and, hence, asymptotically is no more than $n^{1/4}$ (which is much less than the length of any other corridor).

~

Now, denote all partial sums $\sum\limits_{r=0}^{n} f(X_r)$ over first $k$ corridors by $\eta_j, \, 1\le j\le k$. In particular, 
\[
\eta_1 = \sum^{c-1}_{r=0} f(X_r), \quad \eta_2 = \sum^{2c+w-1}_{r=c+w} f(X_r), \;\; \mbox{etc.} 
\]
Note that 
\[
\frac1{\sqrt{n}}|\sum_{r=0}^{n}f(X_r) - \sum_{j=1}^k \eta_j| 
\le C_f \frac{(wk + k)}{\sqrt{n}} \sim  C_f \frac{n^{9/20}+ n^{1/4}}{\sqrt{n}} \to 0, \quad n\to\infty,
\]
uniformly in $\omega \in \Omega$. Hence, it suffices to show that 
\[
\frac1{\sqrt{n}}\sum_{j=1}^k \eta_j \Longrightarrow \eta' \sim {\cal N}(0, \sigma^2).
\]

Note that
\begin{eqnarray*}
n^{-1}{\mathbb E}_{inv} \eta_1^2 
\sim \frac{c}n \sigma^2, \quad n\to\infty, 
\end{eqnarray*}
or, 
\begin{eqnarray}\label{varn}
n^{-3/4}{\mathbb E}_{inv} \eta_1^2 \to  \sigma^2, \quad n\to\infty, 
\end{eqnarray}
and the latter convergence is irrespectively of whether $\sigma^2>0$, or $\sigma^2=0$.

Now, to show the desired weak convergence, let us check the behaviour of the characteristic functions. Due to the Corollary \ref{cor6}, we estimate for any $\lambda \in \mathbb R$, 
\[
|{\mathbb E} (\exp(i\lambda \eta_j) | {\cal F}^{X}_{(j-1)[n/k]}) 
-  {\mathbb E}_{inv} \exp(i\lambda \eta_j)| \le C (1-\kappa)^{[n^{1/5}]}.
\]
So,  by induction, 
\begin{eqnarray} 
&\displaystyle {\mathbb E} \exp(i\frac{\lambda}{\sqrt{n}} \sum_{j=1}^{k+1}\eta_j) 
=  {\mathbb E} \exp(i\frac{\lambda}{\sqrt{n}} \sum_{j=1}^{k}\eta_j) 
{\mathbb E} (\exp(i\frac{\lambda}{\sqrt{n}} \eta_{k+1}) | {\cal F}^{X}_{k[n/k]})
  \nonumber \\ \nonumber \\
&\displaystyle =   {\mathbb E} \exp(i\frac{\lambda}{\sqrt{n}} \sum_{j=1}^{k}\eta_j) 
({\mathbb E}_{inv} (\exp(i\frac{\lambda}{\sqrt{n}} \eta_{k+1}) + O((1-\kappa)^{n^{1/5}})) ) = \ldots
  \nonumber \\ \nonumber \\
&\displaystyle  = {\mathbb E}_{inv} (\exp(i\frac{\lambda}{\sqrt{n}} \eta_{k+1} + O((1-\kappa)^{n^{1/5}}))) ({\mathbb E}_{inv} (\exp(i\frac{\lambda}{\sqrt{n}} \eta_1))^k 
+ O(k(1-\kappa)^{n^{1/5}})). \label{remainder}
\end{eqnarray}
Here $O(k(1-\kappa)^{n^{1/5}})$ is, generally speaking, random and it is a function of $X_{k[n/k]}$, but the module of this random variable does not exceed a nonrandom constant  multiplied by $k(1-\kappa)^{n^{1/5}}$.
We replaced $[n^{1/5}]$ by $n^{1/5}$, which does not change the asymptotic (in)equality. Note that 
\[
O(k(1-\kappa)^{n^{1/5}})) = O(n^{3/4}(1-\kappa)^{n^{1/5}})) \to 0, \quad n\to \infty.
\]

~

Now the idea is to use Taylor's expansion  
\begin{equation}\label{rem2}
{\mathbb E}_{inv} \exp(i\frac{\lambda}{\sqrt{n}} \eta_1) 
= 1 - \frac{\lambda^2}{2n}n^{3/4} \sigma^2 + R_n 
= 1 - \frac{\lambda^2}{2n^{1/4}} \sigma^2 + R_n.
\end{equation}
Here, to prove the desired statement it suffices to estimate accurately the remainder term $R_n$, that is, to show that 
\(R_n = o(n^{-1/4}), \, n\to\infty\).

Since we, actually, transferred the problem to studying an array scheme (as $\eta_1$ itself changes with $n$), we have to inspect carefully this remainder $R_n$.  Due to the Taylor expansion 
we have, 
\begin{eqnarray*}
\mbox{Re}\,\phi(\frac{\lambda}{\sqrt{n}}) = {\mathbb E}_{inv} \cos (\frac{\lambda  \eta_1}{\sqrt{n}}) 
= 1 
-\frac{\lambda^2}{2 n}\mathbb E_{inv}\eta_1^2 + 
\frac{\hat \lambda^3}{6\sqrt{n^3}} {\mathbb E}_{inv}\eta_1^3\sin(\hat\lambda \eta_1),  
\end{eqnarray*}
with some $\hat\lambda$ between $0$ and $\lambda$, 
 and similarly, with some $\tilde\lambda$ between $0$ and $\lambda$,
\begin{eqnarray*}
\mbox{Im}\,\phi(\frac{\lambda}{\sqrt{n}}) = {\mathbb E}_{inv} \sin (\frac{\lambda \eta_1}{\sqrt{n}}) = -
\frac{\tilde \lambda^3}{6\sqrt{n^3}} {\mathbb E}_{inv} \eta_1^3\cos(\tilde\lambda\eta_1). 
\end{eqnarray*}
Here in general $\hat\lambda$ and $\tilde\lambda$ may differ. However, this is not important in our calculus because in any case $|\tilde \lambda| \le |\lambda|$ and $|\hat \lambda| \le |\lambda|$. All we need to do now is to justify a bound 
\begin{equation}\label{just}
|{\mathbb E}_{inv}\eta_1^3| \le K c, 
\end{equation}
with some non-random constant $K$. 
This is a rather standard estimation and we show the details only for completeness. (See similar in \cite{FW, IbragimovLinnik, Khasminsky}, et al.)  It suffices to consider the case $C_f \le 1$, which restriction we assume without loss of generality. 

~

{\bf a.} Consider the case $\mathbb E f(X_k)^3$. We have, clearly, 
\[
|\sum_{k=1}^{c} \mathbb E f(X_k)^3 | \le c.  
\]

~

{\bf b.} For simplicity, denote $f_k = f(X_k)$ and consider the case 
$\mathbb E f_j f_k f_\ell, \, \ell > k>j$. We have, 
\begin{eqnarray*}
&\displaystyle \sum_{j=0}^{c-2} \sum_{k=j+1}^{c-1} \sum_{\ell=k+1}^{c}\mathbb E f_j f_k f_\ell
= \sum_{j=0}^{c-2} \sum_{k=j+1}^{c-1} \sum_{\ell=k+1}^{c}\mathbb E f_j f_k \mathbb E(f_\ell | X_k)
 \\\\
&\displaystyle = \sum_{j=0}^{c-2} \sum_{k=j+1}^{c-1} \sum_{\ell=k+1}^{c}\mathbb E f_j f_k \psi_{k,\ell} q^{\ell-k} \quad \mbox{(here $\psi_{k,\ell} \in {\cal F}^X_{(k)}\equiv \sigma(X_k)$ and $|\psi_{k,\ell}|\le 1$)}.
\end{eqnarray*}
Note that, with a $0\le q < 1$, the expression 
\[
\zeta_{k}:=\sum_{\ell = k+1}^{c-1}  \psi_{k,\ell} q^{\ell-k}
\]
is a random variable, which module is bounded by the absolute constant $(1-q)^{-1}$ and which is ${\cal F}^X_{(k)}$-measurable, i.e., it may be represented as some Borel function of $X_k$. So, we continue the calculus, 
\begin{eqnarray*}
& \displaystyle \sum_{j=0}^{c-2} \sum_{k=j+1}^{c-1} \sum_{\ell=k+1}^{c}\mathbb E f_j f_k f_\ell
= \sum_{j=0}^{c-2} \sum_{k=j+1}^{c-1} \sum_{\ell=k+1}^{c}\mathbb E f_j f_k \zeta_{k} 
= \sum_{j=0}^{c-2} \sum_{k=j+1}^{c-1} \sum_{\ell=k+1}^{c}\mathbb E f_j E(f_k \zeta_{k} | X_j) 
 \\\\
& \displaystyle=  \sum_{j=0}^{c-2} \sum_{k=j+1}^{c-1} 
\mathbb E f_j (\mathbb E f_k \zeta_{k} + \zeta'_{k,j}q^{k-j}) 
= \sum_{j=0}^{c-2} \sum_{k=j+1}^{c-1} \mathbb E f_j  \sum_{k=j+1}^{c-1}  \zeta'_{k,j}q^{k-j} = \sum_{j=0}^{c-2}  \mathbb E f_j  \sum_{k=j+1}^{c-1}  \zeta'_{k,j}q^{k-j},
\end{eqnarray*}
due to $\mathbb E f_j \mathbb E f_k \zeta_{k} = 0$, since $\mathbb E f_j=0$. Here $\zeta'_{k,j}$, in turn, for each $k$ does not exceed by modulus the value $(1-q)^{-1}$ and is ${\cal F}^X_{(j)}$-measurable. Therefore, the inner sum in the last expression satisfies, 
\[
|\sum_{k=j+1}^{c-1}   \zeta'_{k,j}q^{k-j}| 
\le \sum_{k=j+1}^{c-1}   |\zeta'_{k,j}| q^{k-j} 
\le (1-q)^{-1} \sum_{k=j+1}^{c-1}    q^{k-j} 
\le (1-q)^{-2}. 
\]
Thus, 
\[
|\sum_{j=0}^{c-2} \sum_{k=j+1}^{c-1} \sum_{\ell=k+1}^{c}\mathbb E f_j f_k f_\ell| 
\le (1-q)^{-2} \sum_{j=0}^{c-2}  \mathbb E |f_j| 
\le c (1-q)^{-2}, 
\]
as required.

~

{\bf c.} Consider the terms with $\mathbb E f(X_k)^2 f(X_\ell), \, \ell > k$. We estimate, with some (random) $|\psi'_{\ell,k}|\le 1$ and $0\le q<1$, 
\begin{eqnarray*}
|\sum_{k<\ell}^{c-1} \mathbb E f(X_k)^2 \mathbb E(f(X_\ell) | X_k)| 
= |\sum_{k<\ell < c}^{}\mathbb E f(X_k)^2 \psi'_{\ell, k} q^{\ell-k}| \le \frac{c}{1-q}. 
\end{eqnarray*}

~

{\bf d.} Consider the case  $\mathbb E f(X_k)^2 f(X_\ell), \, \ell < k$. We have similarly, for $\ell < k$,
\[
 \mathbb E f^2_k f_\ell = \mathbb E f_\ell \mathbb E(f_k^2 | X_\ell) = \mathbb E f_\ell (\mathbb E f_k^2 + \psi_{\ell, k}'' q^{k-\ell}), \quad |\psi_{\ell, k}|\le 1,
\]
with some (random) $|\psi_{\ell, k}''|\le 1$. 
So, again, 
\begin{eqnarray*}
|\sum_{\ell<k}^{c-1} \mathbb E f_\ell \mathbb E(f_k^2 | X_\ell)| 
= |\sum_{\ell<k<c}^{}\mathbb E f_k^2 \psi''_{\ell, k} q^{k-\ell}| \le \frac{c}{1-q}. 
\end{eqnarray*}

~

{\bf e.} Finally, collecting all intermediate bounds we obtain the bound (\ref{just}), as required:  
\[
|\mathbb E\eta_1^3| \le K c. 
\]
This implies the estimate for the remainder term  $R_n$ in (\ref{rem2}) of the form 
\[
|R_n| \le \frac{c}{n^{3/2}} \sim n^{3/4 - 3/2} = n^{-3/4} = o(n^{-1/4}), 
\]
as required. 
The last detail is to consider the term  ${\mathbb E}_{inv} \exp(i\frac{\lambda}{\sqrt{n}} \eta_{k+1})$ in (\ref{remainder}), for which we have $\sigma^2_{k+1}:= {\mathbb E} \eta^2_{k+1}$ satisfying 
\begin{eqnarray}\label{varnk}
{\mathbb E}_{inv} \eta_{k+1}^2 
= O(n^{3/4} \sigma^2), \quad n\to\infty. 
\end{eqnarray}
This term may be tackled similarly to all others, and, in any case, we get the estimate 
\[
{\mathbb E}_{inv} \exp(i\frac{\lambda}{\sqrt{n}} \eta_{k+1}) 
= 1 + o(1), \quad n\to\infty.
\]
Hence, we eventually get (recall that $c \sim n^{3/4}$), 
\[
\mathbb E\exp(i\frac{\lambda}{\sqrt{n}}\sum_{r=0}^{n} f(X_r)) = 
(1 - \frac{\lambda^2 \sigma^2 c}{2n} + O(\frac{c}{n^{3/2}}))^{n^{1/4}} 
= (1 - \frac{\lambda^2 \sigma^2}{2n^{1/4}} + O(\frac{1}{n^{3/4}}))^{n^{1/4}} 
\to \exp(-\lambda^2 \sigma^2/2),
\]
which is the characteristic function for the Gaussian distribution ${\cal N}(0,\sigma^2)$, as required. 

~

{\bf II.} The case $\sigma^2 = 0$ is considered absolutely similarly. 
Namely, with a practically identical arguments we get, now with $\sigma^2=0$, 
\[
\mathbb E\exp(i\frac{\lambda}{\sqrt{n}}\sum_{r=0}^{n} f(X_r)) = 
(1 - \frac{\lambda^2 \sigma^2 c}{2n} + O(\frac{c}{n^{3/2}}))^{n^{1/4}} 
=  (1  + O(\frac{1}{n^{3/4}}))^{n^{1/4}} 
\to 1,\]
which is the characteristic function for the degenerate Gaussian distribution ${\cal N}(0,0)$, as required.
Hence, the Theorem \ref{thm_clt1} is proved.

\section{Coupling method for Markov chain: simple version}\label{sec:coupling_s}
Concerning coupling method, it is difficult to say who exactly invented this method. The common view -- shared by the author of these lecture notes -- is that it was introduced by W. Doeblin \cite{Doeblin}, even though he himself refers to some ideas of Kolmogorov with relation to the study of ergodic properties of Markov chains. Leaving this subject to the historians of mathematics, let us just mention that there are quite a few articles and monographs 
where this method is presented \cite{Griffeath, Lindvall, Nummelin, Thorisson}, et al. Also there many papers and books where this or close method is used for further investigations without being explicitly named, see, e.g., \cite{Bor98}. This method itself provides ``another way'' to establish geometric convergence in the Ergodic theorem. In the simple form as  in this section this method has limited applications; however, in a more elaborated version -- see the section \ref{sec:gcm} below -- it is most useful, and applicable to a large variety of Markov processes including  rather general diffusions, providing not necessarily geometric rates of convergence but also much weaker rates in non-compact spaces.

By simple coupling for two random variables $X^1, X^2$ we understand the situation where both $X^1, X^2$ are defined on the same probability space and 
\[
\mathbb P(X^1 = X^2)>0.
\]
\noindent
Consider a Markov chain $(X_n, \, n=0, 1, \ldots)$.  
In fact, this simple ``Doeblin's'' version of coupling provides bounds of convergence which are far from optimal in most cases. (By ``far from optimal'' we understand that the constant under the exponential is too rough.) Yet, its advantage is its simplicity and, in particular, no change of the initial probability space. From the beginning we need two ``independent'' probability spaces, $(\Omega^1, {\cal F}^1, \mathbb P^1)$, and $(\Omega^2, {\cal F}^2, \mathbb P^2)$, and the whole construction runs on the direct product of those two: 
\[
(\Omega, {\cal F}, \mathbb P):= (\Omega^1, {\cal F}^1, \mathbb P^1)  \times (\Omega^2, {\cal F}^2, \mathbb P^2).
\]
This space \((\Omega, {\cal F}, \mathbb P)\) will remain unchanged in this section. We assume that there are two Markov processes $(X^1_n)$ on \((\Omega^1, {\cal F}^1, \mathbb P^1)\) and $(X^2_n)$ on \((\Omega^2, {\cal F}^2, \mathbb P^2)\), correspondingly, with the same transition probability matrix ${\cal P}=(p_{ij})_{i,j\in  S}$ satisfying the ``simple ergodic assumption'', 
\begin{equation}\label{sea}
\kappa_0:= \min_{i,j}p_{ij} > 0.
\end{equation}
Naturally, both processes are defined on \((\Omega, {\cal F}, \mathbb P)\) as follows, 
\[
X^1_n(\omega) = X^1_n(\omega^1, \omega^2) := X^1_n(\omega^1), \quad \&  \quad X^2_n(\omega) = X^2_n(\omega^1, \omega^2) := X^2_n(\omega^2). 
\]
We will need some (well-known) auxiliary results. Recall that given a filtration $({\cal F}_n)$, stopping time is any random variable $\tau < \infty$ a.s. with values in ${\mathbb Z}_+$ such that for any $n\in {\mathbb Z}_+$, 
\[
(\omega: \, \tau > n) \in {\cal F}_n.
\]

\noindent
In most of textbooks on Markov chains the following Lemma may be found (see, e.g., \cite{Wentzell}). 
\begin{lemma}\label{smp}
Any Markov chain (i.e., a Markov process with discrete time) is strong Markov. 
\end{lemma}
Consider a new process composed from two, $X_ n := (X^1_n, X^2_n)$, evidently, with two independent coordinates. Due to this independence, the following Lemma holds true. 
\begin{lemma}\label{xxx}
The (vector-valued) process $(X_n)$ is a (homogeneous) Markov chain; hence, this chain is also strong Markov. 
\end{lemma}

In the following main result of this section, $\mu$ stands for the (unique) stationary distribution of our Markov chain $(X^1_n)$ (as well as of $(X^2_n)$). 
\begin{theorem}\label{weak_kappa}
For any initial distribution $\mu_0$,
 
\begin{equation}\label{weakexp}
\sup_{A} |P_{\mu}(n,A) - \mu(A)| \le (1-\kappa_0)^n.
\end{equation}
\end{theorem}
Let us emphasize again that the bound may be not optimal; however, the advantage is that the construction of coupling here does not require any change of the probability space.

~

\noindent
{\em Proof of Theorem \ref{weak_kappa}.}
Recall that a new Markov chain $X_ n := (X^1_n, X^2_n)$ with two independent coordinates is strong Markov. Let
\[
\tau:= \inf(n\ge 0: \; X^1_n = X^2_n). 
\]
We have seen that $\mathbb P(\tau<\infty)=1$. More than that, from Markov property it follows for any $n$ by induction (with a random variable called {\em indicator,} $1(A)(\omega) = 1$ if $\omega \in A$ and $1(A)(\omega) = 0$ if $\omega \not\in A$), 
\begin{eqnarray}\label{bd0}
&\displaystyle {\mathbb P}(\tau>n) = \mathbb E1(\tau>n) = \mathbb E\prod_{k=1}^{n}1(\tau>k)
 \nonumber \\ \nonumber\\ \nonumber 
&\displaystyle = \mathbb E \left(\mathbb  E(\prod_{k=1}^{n}1(\tau>k) | {\cal F}_{n-1})\right) 
= \mathbb E \left(\prod_{k=1}^{n-1}1(\tau>k) \mathbb E(1(\tau>n) | {\cal F}_{n-1})\right)
 \\\nonumber\\
&\displaystyle \le \mathbb E\prod_{k=1}^{n-1}1(\tau>k) (1-\kappa_0) =  (1-\kappa_0) \mathbb E\prod_{k=1}^{n-1}1(\tau>k) \le \; \stackrel{\mbox{\small (induction)}}{\ldots} \; \le (1-\kappa_0)^n.
\end{eqnarray}
Define
\[
X^3_n := X^1_n 1(n < \tau) + X^2_n 1(n\ge \tau). 
\]
Due to the strong Markov property, $(X^3)$ is also a {\em Markov chain} and it is  {\em equivalent} to $(X^1)$ -- that is, they both have the same distribution in the space of trajectories.  
This follows from the fact that at $\tau$ which is a stopping time the processes follows $X^2$, so that it uses the same transition probabilities while choosing the next state at $\tau+1$ and further.

~

Now, here is the most standard and most frequent calculus in most of works on coupling method, or where this method is used (recall that all the processes $X^1, \, X^2, \, X^3$ are defined on the same probability space): for any $A\in {\cal S}$, 
\begin{eqnarray*}
&\displaystyle |{\mathbb P}(X^1_n\in A) - {\mathbb P}(X^2_n\in A)| = |{\mathbb P}(X^3_n\in A) - {\mathbb P}(X^2_n\in A)|
 \\\\
&\displaystyle = |\mathbb E1(X^1_n\in A) - \mathbb E1(X^2_n \in A)| 
= |\mathbb E(1(X^3_n \in A) - 1(X^2_n\in A))| 
 \\\\
&\displaystyle =   |\mathbb E(1(X^3_n\in A) - 1(X^2_n\in A)) 1(\tau>n) + \mathbb E(1(X^3_n\in A) - 1(X^2_n\in A)) 1(\tau \le n)|
 \\\\
&\displaystyle \stackrel{(*)}{=}  |\mathbb E(1(X^3_n\in A) - 1(X^2_n\in A)) 1(\tau>n)| \le  |\mathbb E(1(X^3_n\in A) - 1(X^2_n\in A)) 1(\tau>n)| 
 \\\\
&\displaystyle \le \mathbb E|1(X^3_n\in A) - 1(X^2_n\in A)| 1(\tau>n) 
\le \mathbb E 1(\tau>n) = {\mathbb P}(\tau>n) \stackrel{(\ref{bd0})}{\le} (1-\kappa_0)^n. 
\end{eqnarray*}
Note that the final bound is uniform in $A$. Here the equality (*) is due to the trivial fact that since $n\ge \tau$, the values of $X^3_n$ and $X^2_n$ coincide, so either \(1(X^3_n\in A) = 1(X^2_n\in A) = 0\), or  \(1(X^3_n\in A) = 1(X^2_n\in A) = 1\) simultaneously on each $\omega$, which immediately implies that \((1(X^3_n\in A) - 1(X^2_n\in A)) 1(\tau \le n) = 0\). So, the Theorem \ref{weak_kappa} is proved.

\section{A bit of large deviations}\label{sec:ld}
In this section assume 
$$
{\mathbb E}_{inv} f(X_0) = 0. 
$$
We will be interested in the existence and properties of the limit, 
\begin{equation}\label{hb}
\lim_{n\to\infty} \frac1{n} \ln {\mathbb E}_{x} \exp(\beta \sum_{k=0}^{n-1} f(X_k)) =: H(\beta). 
\end{equation}
Note that we do not use $x$ in the right hand side because in ``good cases'' -- as below -- the limit does not depend on the initial state. 
Denote 
$$
H_n(\beta,x) := \frac1{n} \ln {\mathbb E}_{x} \exp(\beta \sum_{k=0}^{n-1} f(X_k)), 
$$
and define the operator $T = T^\beta$ acting on functions on $S$ as follows, 
\[
T^\beta h(x) := \exp(\beta f(x)) {\mathbb E}_x h(X_1),  
\]
for any function $h$ defined on $S$. 
Note that 
\[
{\mathbb E}_{x} \exp(\beta \sum_{k=0}^{n-1} f(X_k)) =  (T^\beta)^n h(x), 
\]
with $h(x) \equiv 1$. 
Indeed, for $n=1$ this coincides with the defiition of $T^\beta$. 
Further, for $n >1$ due to the Markov property by induction,
\begin{eqnarray*}
&\displaystyle {\mathbb E}_{x} \exp(\beta \sum_{k=0}^{n-1} f(X_k)) = {\mathbb E}_{x} {\mathbb E}_{x} (\exp(\beta \sum_{k=0}^{n-1} f(X_k))| X_{1}) 
 \\\\
&\displaystyle = \exp(\beta f(x)) {\mathbb E}_{x}  
{\mathbb E}_{x}(\exp(\beta \sum_{k=1}^{n-1} f(X_k)) | X_{1})
 \\\\
&\displaystyle =  \exp(\beta f(x)) {\mathbb E}_{x}  (T^\beta)^n h(x) = T^\beta (T^\beta)^{n-1} h(x) = (T^\beta)^{n} h(x),
\end{eqnarray*}
as required. 
So, the function $H_n$ can be rewritten as
$$
H_n(\beta,x) = \frac1{n} \ln (T^\beta)^nh(x),  
$$
($h(x) \equiv 1$).
It is an easy exercise to check that the function $H_n(\beta,x)$ is convex in $\beta$. Hence, if the limit exists, then the limiting function $H$ is also convex. Now recall the following classical and basic result about positive matrices. 
\begin{theorem}[Perron--Frobenius]\label{thm-pf}
Any positive quadratic matrix (i.e., with all entries positive) has  a positive eigenvalue $r$ called its spectral radius, which is strictly greater than the moduli of the rest of the spectrum, this eigenvalue is simple, and its corresponding eigenfunction (eigenvector) has all positive  coordinates. 
\end{theorem}
In fact, this result under the specified conditions is due to Perron, while Frobenius extended it to non-negative matrices. We do not discuss the details of this difference and how it can be used. Various presentations may be found, in particular, in \cite{Karlin, KLS, Seneta1}. 
As an easy corollary, the Theorem \ref{thm-pf} implies the existence of the limit in (\ref{hb}) -- which, as was promised, does not depend on $x$ -- with, 
\begin{equation}\label{hr}
H(\beta) = \ln r(\beta), 
\end{equation}
where $r(\beta)$ is the spectral radius of the operator $T^\beta$,
see, for example, \cite[Theorem 7.4.2]{FW}. (Emphasize that in the proof of this theorem it is important that the eigenvector corresponding to the spectral radius is strictly positive, i.e., it has all positive components.) More than that, in our case it follows from the theorem about analytic properties of simple eigenvalues $r(\beta)$ is analytic, see, e.g., \cite{Kato}. Therefore, $H(\beta)$ is analytic, too. Also, clearly, analytic is $H_n$ as a function of the variable $\beta$. Then it follows from the properties of analytic (or convex) functions that convergence $H_n(\beta,x) \to H(\beta)$ implies that also
\[
H'_n(0,x) \to H'(0), \quad n\to\infty,  
\]
where by $H'_n$ we understand the derivative with respect to $\beta$.
On the other hand, we have, 
\begin{eqnarray*}
& \displaystyle H_n'(0,x)  = \frac{\partial}{\partial \beta}
\left(\frac1{n} \ln {\mathbb E}_{x} \exp(\beta \sum_{k=0}^{n-1} f(X_k))\right)|_{\beta = 0}
 \\\\
& \displaystyle = \frac1{n} \, \frac{{\mathbb E}_{x} \left(\sum\limits_{k=0}^{n-1} f(X_k) \exp(\beta \sum\limits_{k=0}^{n-1} f(X_k))\right)} {{\mathbb E}_{x}  \exp(\beta \sum_{k=0}^{n-1} f(X_k))}|_{\beta = 0}
= \frac1n {\mathbb E}_{x}  \sum\limits_{k=0}^{n-1} f(X_k). 
\end{eqnarray*}
So, due to the Law of Large Numbers, 
$$
H_n'(0,x)  =\frac1n {\mathbb E}_{x}  \sum_{k=0}^{n-1} f(X_k) \to {\mathbb E}_{inv} f(X_0) = 0.
$$
Hence, 
$$
H'(0) = {\mathbb E}_{inv} f(X_0) = 0.  
$$
Also, again due to the analyticity, 
\[
H''_n(0,x) \to H''(0), \quad n\to\infty. 
\]
On the other hand, due to (\ref{sig}) and (\ref{varn}), 
$$
H''_n(0)  =\frac1n {\mathbb E}_{x}  (\sum_{k=0}^{n-1} f(X_k))^2 \to \sigma^2, \quad n\to\infty. 
$$
Hence, 
$$
H''(0) = \sigma^2.
$$
This last assertion will not be used in the sequel. 

~

\noindent
Let us state it all as a lemma. 
\begin{lemma}\label{lem-hp}
There exists a limit $H(\beta)$ in (\ref{hb}). This function $H$ is convex and differentiable, and, in particular, 
$$
H'(0) = 0, \quad H''(0) = \sigma^2. 
$$
\end{lemma}

Actually, we will not use {\em large deviations} (LDs) in these lecture notes, except  for the Lemma \ref{lem-hp}, which is often regarded as a preliminary auxiliary result in large deviation theory. Yet, once the title of the section uses this term, let us state one simple inequality of LD  type. 
Recall that ${\mathbb E}_{inv}f(X_0) = 0$ in this section. 

\begin{proposition}\label{ld_pro}
Let
\begin{equation}\label{ltr}
L(\alpha) := \sup_{\beta}(\alpha \beta - H(\beta)), \quad \tilde L(\alpha):= 
\limsup_{\delta \downarrow 0}
(L(\alpha - \delta)1(\alpha>0) + L(\alpha + \delta)1(\alpha<0)).
\end{equation}
Then under the assumptions of the Ergodic Theorem  \ref{thm_erg1} for any $\epsilon > 0$, 
\begin{equation}\label{ld1}
\limsup_{n\to\infty} \frac1n \ln \mathbb P_x(\frac1n \sum_{k=0}^{n-1} f(X_k)  \ge \epsilon) \le - \tilde L(\epsilon).
\end{equation}
\end{proposition}
The function $L$ is called Legendre transformation of the function $H$. It is convex and lower semicontinuous; see \cite{Rocka} about this and more general transformations (e.g., where $H$ is convex but not necessarily smooth -- in which case $L$ is called Fenchel--Legendre's transformation). Notice that ``usually'' (everywhere where $L(\alpha)<\infty$) in (\ref{ld1})  there is a limit instead of $\limsup$, and this limit equals the right hand side, and both $\tilde L(\epsilon) = L(\epsilon)>0$; the latter is certainly true, at least, for small $\epsilon >0$ if $\sigma^2 >0$. Also, for positive $\alpha$, of course, it suffices to use the definition $\displaystyle \tilde L(\alpha):= \limsup_{\delta \downarrow 0}
L(\alpha - \delta)$ in (\ref{ltr}). However, this simple result does not pretend to be even an introduction to large deviations, about which theory see \cite{BorMog, DZ, FK, FW, GuVer93, Puh, Var}, et al. In the next sections the Proposition \ref{ld_pro} will not be used: all we will need is the limit in (\ref{hb}) due to the Lemma \ref{lem-hp} and some its properties which will be specified. 

~

\noindent
{\em Proof of Proposition \ref{ld_pro}.} We have for any $0<\delta < \epsilon$, by Chebyshev--Markov's exponential inequality with any $\lambda>0$, 
\begin{eqnarray*}
\mathbb P_x(\frac1n \sum_{k=0}^{n-1} f(X_k)  \ge \epsilon) = \mathbb P_x(\exp(\lambda\sum_{k=0}^{n-1} f(X_k))\ge \exp( n\lambda \epsilon))
 \\\\
\le \exp(-n\lambda \epsilon) 
\mathbb E_x\exp(\lambda\sum_{k=0}^{n-1} f(X_k)) \stackrel{(\ref{hb})}{\le} \exp(-n(\lambda (\epsilon -\delta) +H(\lambda))),  
\end{eqnarray*}
if $n$ is large enough. The first and the last terms here with the inequality between them can be rewritten equivalently as 
\[
\frac1n \ln \mathbb P_x(\frac1n \sum_{k=0}^{n-1} f(X_k)  \ge \epsilon) \le -\lambda (\epsilon -\delta) +H(\lambda), 
\]
for $n$ large enough. So, we have, 
\[
\limsup\limits_{n\to\infty} \frac1n \ln \mathbb P_x(\frac1n \sum_{k=0}^{n-1} f(X_k)  \ge \epsilon) \le -(\lambda (\epsilon -\delta) - H(\lambda)). 
\]
Since this is true for any $\lambda >0$, we also get, 
\begin{eqnarray*}
\limsup\limits_{n\to\infty} \frac1n \ln \mathbb P_x(\frac1n \sum_{k=0}^{n-1} f(X_k)  \ge \epsilon) 
\le -\sup_{\lambda>0} (\lambda (\epsilon -\delta) -H(\lambda)),  
\end{eqnarray*}
However, since $H(0)=0$ and $H'(0) = 0$ and due to the convexity of $H$, the supremum on {\em all $\lambda\in \mathbb R$} here on positive $\epsilon - \delta$ is attained at $\lambda>0$, i.e., 
$$
\sup_{\lambda>0}(\lambda (\epsilon -\delta) -H(\lambda)) = \sup_{\lambda\in \mathbb R}(\lambda (\epsilon -\delta) -H(\lambda)) \equiv L(\epsilon - \delta).
$$ 
Thus, the left hand side in (\ref{ld1}) does not exceed the value $- \limsup\limits_{\delta \downarrow 0}L(\epsilon - \delta) = -\tilde L(\epsilon)$, as required. The Proposition \ref{ld_pro} is proved.

\section{Dynkin's formulae}
Let $L$ be a generator of some Markov chain on a finite state space $S$, that is, for any function $u$ on $S$, 
\begin{equation}\label{gen}
Lu : = \mathbb E_x u(X_1) - u(x) \equiv {\cal P}u(x) - u(x). 
\end{equation}
Recall that here ${\cal P}$ is the transition probability matrix of the corresponding Markov chain $(X_n)$, a function $u$ on $S$ is considered as a column-vector, ${\cal P} u$ is this matrix multiplied by this vector, and ${\cal P}u(x)$ is the $x$-component of the resulting vector.
Note that such difference operators are discrete analogues of elliptic differential operators of the second order studied extensively, in particular, in mathematical physics. What makes them the analogues is that both are generators of Markov processes, either in discrete or in continuous time; also, it may be argued about limiting procedures approximating continuous time processes by discrete ones. Yet, the level of this comparison here is, of course, intuitive and we will not try to justify in any way, or to explain it further. 

As usual in these lecture notes, we will assume that the corresponding process $(X_n)$ satisfies the Ergodic Theorem \ref{thm_erg1}. 
The Poisson equation for the operator $L$ from (\ref{gen}) is as follows: 
\begin{equation}\label{pe}
Lu(x)  = -f(x), \quad x\in S. 
\end{equation}
This equation may be studied {\em with} or {\em without} some boundary and certain boundary conditions. 
The goal of this chapter is to present how such equations may be solved probabilistically. This simple study may be also considered as an introduction to the  Poisson equations  for elliptic  differential operators. 
We start with {\bf Dynkin's formula} or Dynkin's identity. 
\begin{theorem}[Dynkin's formula 1]\label{Df1}
On the finite state space $S$, for any function $h$ and any $n = 1,2,\ldots$, 
\begin{equation}\label{de}
\mathbb E_x h(X_n) = h(x) + \sum_{k=0}^{n-1} \mathbb E_x Lh(X_k), \quad n\ge 0.  
\end{equation}

\end{theorem}

\noindent
{\em Proof.} {\bf 1.} For $n=1$ the formula (\ref{de}) reads, 
$$
\mathbb E_x h(X_1) = h(x) + 
Lh(x), 
$$
where $x$ is a non-random initial value of the process. 
Hence, by inspection, the desired identity for $n=1$ is equivalent to the definition of the generator in (\ref{gen}). 

~

\noindent
{\bf 2.} For the general case $n$, the desired formula follows by induction. Indeed, assume that the formula (\ref{de}) holds true for some $n=k$ and check it for $n=k+1$. We have, 
\begin{eqnarray*}
&\displaystyle \mathbb E_x h(X_{n+1}) 
=  \mathbb E_x h(X_{n+1}) - \mathbb E_x h(X_{n}) + \mathbb E_x h(X_{n})  
 \\\\
&\displaystyle =  \mathbb E_x \mathbb E_x (h(X_{n+1}) - h(X_{n}) | X_n)+ \mathbb E_x h(X_{n})  
 \\\\
&\displaystyle =  \mathbb E_x \mathbb E_x (Lh(X_n)|X_n)+   h(x) + \sum_{k=0}^{n-1} \mathbb E_x Lh(X_k)
 \\\\
&\displaystyle  =  \mathbb E_x Lh(X_n) +   h(x) + \sum_{k=0}^{n-1} \mathbb E_x Lh(X_k)
= h(x) + \sum_{k=0}^{n} \mathbb E_x Lh(X_k). 
\end{eqnarray*}
So, the formula (\ref{de}) for all values of $n$ follows by induction. The Theorem \ref{Df1} is proved. 

\section{Stopping times and martingales: reminder}
\begin{definition}\label{Def-mart}
Filtration $({\cal F}_n, \, n=0,1,\ldots)$ is a family of increasing sigma-fields on a probability space $(\Omega, {\cal F}, \mathbb P)$ completed with respect to the measure $\mathbb P$ (that is, each ${\cal F}_n$ contains each subset of all $\mathbb P$-zero sets from ${\cal F}$). 
The process $(M_n)$ is called a martingale with respect to a filtration $({\cal F}_n)$ iff $\mathbb E M_n < \infty$ and $\mathbb E(M_{n+1} |{\cal F}_n) = M_n$ (a.s.). 
\end{definition}

\begin{definition}\label{stime}
A random variable $\tau < \infty$ a.s. with values in ${\mathbb Z}_+$ is called a stopping time with respect to a filtration $({\cal F}_n)$ iff for each $n\in {\mathbb Z}_+$ the event $(\tau>n)$ is measurable with respect to ${\cal F}_n$.
\end{definition}

\noindent
It is recommended to read about simple properties of martingales and stopping times in one of the textbooks on stochastic processes, e.g., \cite{Kry-ln}. We will only need the following classical result about stopped martingales given here without proof. 

\begin{theorem}[Doob]
Let $(M_n)$ be a martingale and let $\tau$ be a stopping time with respect to a filtration $({\cal F}_n)$. Then $(\tilde M_n:= M_{n\wedge \tau})$ is also a martingale. 
\end{theorem}
In terms of martingales, the first Dynkin's formula may be re-phrased as follows. 
\begin{theorem}[Dynkin's formula 2]\label{df2}
On the finite state space $S$, for any function $h$ and any $n = 1,2,\ldots$, the process 
\begin{equation}\label{de2}
M_n:= h(X_n) - h(x) - \sum_{k=0}^{n-1} Lh(X_k), \quad n\ge 0, 
\end{equation}
is a martingales with respect to the natural filtration ${\cal F}^X_n$ ``generated'' by the process $X$. Vice versa, if the process $M_n$ from (\ref{de2}) is a martingale then (\ref{de}) holds true.
\end{theorem}
\noindent
{\em Proof.} The inverse statement is trivial. The main part follows 
due to the Markov property, 
\begin{eqnarray*}
{\mathbb E}(M_n | {\cal F}_{n-1}) = {\mathbb E} (h(X_n) | X_{n-1}) - h(x) - \sum_{k=0}^{n-1} Lh(X_k)
 \\\\
= {\cal P}h (X_{n-1}) - Lh(X_{n-1}) 
- h(x) - \sum_{k=0}^{n-2} Lh(X_k)
 \\\\
= h(X_{n-1}) - h(x) - \sum_{k=0}^{n-2} Lh(X_k) = M_{n-1}.  
\end{eqnarray*}
The Theorem \ref{df2} is thus proved.

\begin{lemma}[Dynkin's formula 3]\label{dynkin}
Let $\tau$ be a stopping time with 
$$
\mathbb E_x\tau < \infty, \quad \forall \; x\in S. 
$$
Then for any function $h$ on $S$, 
\[
\mathbb E_x h(X_\tau) =  h(x) + \mathbb E_x\sum_{k=0}^{\tau-1}  Lh (X_k). 
\]
\end{lemma}
\noindent
{\em Proof.} Follows straightforward from the Theorem \ref{df2}.

\section{Poisson equation without a potential}\label{sec:pnp}
\subsection{Introduction}

Here we consider the following discrete Poisson  equation {\it without a potential,}
\begin{equation}\label{eq02}
Lu(x) \equiv {\cal P}u(x) -u(x)  = -f(x).  
\end{equation}
In  the next section a similar discrete equation {\it with a potential $c = c(x), \, x\in S$,} will be studied,  
\begin{equation}\label{eq020}
L^c u(x) := \exp(-c(x)){\cal P}u(x) -u(x)= -f(x),
\end{equation}
firstly because it is natural for PDEs -- and here we present an easier but similar discrete-time theory --  
and secondly with a hope that it may be also 
useful for some further extensions, as it already happened with equations without a potential. 
Let $\mu$ be, as usual,  the (unique) invariant probability measure of the process $(X_n, \, n\ge 0)$.

\subsection{Poisson equation (\ref{eq02}) with a boundary} 

Firstly, we consider Poisson equation with a non-empty boundary,   
\begin{equation}\label{peb}
L u (x) = -f(x), \;\; x\in S \setminus \Gamma, \quad u(x) = g(x), \; x\in \Gamma, 
\end{equation}
where $\Gamma \subset S$, $\Gamma \not = \emptyset$.
If the right hand side equals zero, this equation is called the Laplace equation with Dirichlet boundary conditions:
\begin{equation}\label{leb}
L u (x) = 0, \;\; x\in S \setminus \Gamma, \quad u(x) = g(x), \; x\in \Gamma.
\end{equation}

\noindent
Let 
\[
\tau:= \inf(n\ge 0:\; X_n \in \Gamma), 
\]
and denote
\begin{equation}\label{pe1_sol}
v(x) := \mathbb E_x \left(\sum_{k=0}^{\tau-1}f(X_k) + g(X_\tau)\right).
\end{equation}
Recall that under our assumptions on the process, necessarily ${\mathbb E}_x\tau < \infty$.

~

For the uniqueness, we would need a {\bf maximum principle}, which holds true for the Laplace equation (recall that we always assume $\min\limits_{i,j} p_{ij}>0$): 

\begin{lemma}[Maximum principle]\label{lemp}
If the function $u$ satisfies the equation (\ref{leb}), then the maximal value (as well as minimal) of this function is necessarily attained at the boundary~$\Gamma$. 
\end{lemma}

\noindent
{\em Proof.} Since $Lu(x)=0$ for any $x \not \in \Gamma$, we have 
\begin{equation}\label{le}
u(x) = {\cal P}u(x),
\end{equation}
for such $x$. In other words, the value $u(x)$ is equal to the average of the values $u(y)$ at {\bf all} other $y\in S$ with some positive weights, due to the assumption $\min\limits_{ij}p_{ij}>0$. However, if a maximal value, say, $M$, is attained by  $u$ {\em not} at the boundary, say, $u(x_0)=M$, $x_0 \not\in\Gamma$, and if at least one value on $\Gamma$ (or, actually, anywhere else) is strictly less than $M$, then we  get a contradiction, as the equality $\sum\limits_{y \in  S} p_{xy} v(y) = M$ with all $v(y)\le M$ and with at least, one $v(y)<M$ is impossible. Similar arguments apply to the minimal value of $u$. This proves the Lemma \ref{lemp}.

\begin{theorem}\label{thm_pe1sol}
The function $v(x)$ given by the formula (\ref{pe1_sol}) is a unique solution of the Poisson equation (\ref{peb}). 
\end{theorem}
\noindent
{\em Proof.} {\bf 1.} The boundary condition $v(x) = g(x)$ on $x \in \Gamma$  is trivial because $\tau = 0$ in this case. 

~

\noindent
{\bf 2.} 
Let $x\not\in \Gamma$. Then $\tau \ge 1$. We have, due to the Markov property,
\begin{eqnarray*}
&\displaystyle v(x) = f(x) +  \sum_{y} {\mathbb E}_x 1(X_1=y) {\mathbb E}_{y} \left(\sum_{k=0}^{\tau-1} f(X_k) + g(X_\tau)\right)
 \\\\
&\displaystyle = f(x) +  \sum_{y} p_{xy} v(y) =  f(x) + {\mathbb E}_x v(X_1).
\end{eqnarray*}
From this, it follows clearly the statement about solving the equation,
$$
 Lv(x) = {\mathbb E}_x v(X_1) - v(x) = -f(x).
$$

\noindent
{\bf 3.} Uniqueness follows from the maximum principle. Indeed, let $v^1$ and $v^2$ be two solutions. Then 
$$
u(x):= v^1(x) - v^2(x) = 0, \quad \forall \;\; x\in \Gamma.
$$
Also, at any $x\not\in \Gamma$, 
$$
Lu(x) = Lv^1(x) - Lv^2(x) = 0. 
$$
Hence, by virtue of the Lemma \ref{lemp},  both maximal and minimal values of the function $u$ are attained at the boundary $\Gamma$. However, at the boundary both these values are equal to zero. Therefore, 
$$
u(x)= 0, \quad \forall \; x\in S,  
$$
that is, $v^1-v^2 \equiv 0$, as required. This completes the proof of the Theorem \ref{thm_pe1sol}.

\subsection{Poisson equation (\ref{eq02}) without a boundary}
Consider the equation on the whole $S$, 
\begin{equation}\label{pew}
L u (x) = -f(x), \;\; x\in S.
\end{equation}
We will need an additional assumption on $f$ called ``centering''. This condition is a close analogue of the subtraction in the standardization for a CLT.

\begin{Assumption}[Centering]\label{As4}
It is assumed that the function $f$ satisfies the condition, 
\begin{equation}\label{cent}
{\mathbb E}_{inv}f(X_0) \equiv  \sum_x f(x)\mu(x) = 0, 
\end{equation}
where $\mu$ is the (unique) invariant measure of the process $X$. 
\end{Assumption}

\begin{theorem}\label{pew-sol}
Under the assumption (\ref{cent}), the equation 
(\ref{pew}) has a solution $u$, 
which is unique up to an additive constant. 
This solution is given by the formula 
\begin{equation}\label{pesol1}
u(x) = \sum_{k=0}^{\infty} {\mathbb E}_{x} f(X_k). 
\end{equation}
The solution $u$ from (\ref{pesol1}) itself satisfies the centering condition, 
\begin{equation}\label{centu}
\sum u(x) \mu(x) = 0. 
\end{equation}

\end{theorem}
Note that the ``educated guess'' about a solution represented by the formula (\ref{pesol1}) may be deduced from the comparison with (\ref{pe1_sol}) where, so to say, we want to drop the terminal summand $g$ as there is no boundary and to replace $\tau$ by infinity; naturally, expectation and  summation should be interchanged. 
Also,  in the present setting the idea based on considering the series for $(I- {\cal P})^{-1}$ on centered functions may be applied. Yet, we would like to avoid this way because in a more general ``non-compact'' situation a {\em polynomial} convergence of the series in (\ref{pesol1}) would also suffice, and, hence, this approach looks more general.

~

\noindent
{\em Proof of Theorem \ref{pew-sol}. }
{\bf 1. Convergence.} Follows straightforward from the Corollary \ref{cor6}. This shows that the function $u(x)$ defined in (\ref{pesol1}) is everywhere finite. 
  
~

\noindent
{\bf 2. Satisfying the equation.} From the Markov property,
\begin{eqnarray*}
u(x) = f(x) +  \sum_{y} {\mathbb E}_x 1(X_1=y) {\mathbb E}_{y} \sum_{k=0}^{\infty} f(X_k)
 \\\\
= f(x) +  \sum_{y} p_{xy} v(y) =  f(x) + {\mathbb E}_x u(X_1).
\end{eqnarray*}
From this, it follows clearly the statement,
$$
 Lu(x) = {\mathbb E}_x u(X_1) - u(x) = -f(x).
$$

~

\noindent
{\bf 3. Uniqueness.} Let $u^1$ and $u^2$ be two solutions both satisfying the moderate growth and centering. Denote $v = u^1 - u^2$. Then
\[
Lv = 0.
\]
By virtue of Dynkin's formula  (\ref{de}), 
\[
{\mathbb E}v(X_n) - v(x) = 0. 
\]
However, due to the  Corollary \ref{cor6}, 
\[
{\mathbb E}_x v(X_n) \to {\mathbb E}_{inv} v(X_0) =  0. 
\]
Hence, 
\[
v(x) \equiv 0, 
\]
as required. 

~

\noindent
{\bf 4. Centering.} We have, due to a good convergence -- see the Corollary \ref{cor6} -- and Fubini's theorem, and since measure $\mu$ is stationary, and finally because $f$ is centered, 
\begin{eqnarray*}
&\displaystyle \sum_x u(x)\mu(x)  = \sum_x \mu(x)\sum_{k=0}^{\infty} {\mathbb E}_{x} f(X_k)
 \\\\
&\displaystyle = \sum_{k=0}^{\infty}\sum_x \mu(x) {\mathbb E}_{x} f(X_k) =  \sum_{k=0}^{\infty} {\mathbb E}_{inv} f(X_k) = 0.
\end{eqnarray*}
The Theorem  \ref{pew-sol} is proved.

\section{Poisson equation with a potential}\label{sec:pp}
Let us remind the reader  that the case $|S|<\infty$ is under consideration. 

\subsection{Equation (\ref{eq020})}
Recall the equation (\ref{eq020}), 
\[
\exp(-c(x)){\cal P}u(x) - u(x)= -f(x).
\]
A natural candidate for the solution is the function\\
\begin{equation}\label{pe2-sol}
 u(x) := \sum_{n=0}^\infty {\mathbb E}_x \exp\left(-\sum_{k=0}^{n-1} c(X_k)
\right) f(X_k),
\end{equation}
{\bf provided that this expression is well-defined}. Naturally on our finite state space $S$ both $f$ 
and $c$ bounded. 
Denote
\[
\varphi_n := \sum\limits_{k=0}^{n} c(X_k), \quad \varphi_{-1}=0, 
\]
and
\[
L^c := \exp(-c(x)){\cal P} - I, 
\]
that is,
\[
L^c u(x) := \exp(-c(x)){\cal P}u(x) - u(x).
\]

\smallskip

\noindent
We can tackle several cases, and the most interesting one in our view is where $c(x) = \varepsilon  c_1(x)$, $\varepsilon > 0$ 
\emph{small} and $\displaystyle \bar c_1:= \sum\limits_{x} c_1(x) \mu(x) > 0$.
Denote also $\bar c = \sum\limits_x c(x) \mu(x)$.

\subsection{Further Dynkin's formulae}

\begin{lemma}[Dynkin's formula 4]\label{dynkin4}
\begin{equation}\label{dyn11}
\mathbb E_x \exp(-\varphi_{n-1})\, h(X_n) =  h(x) + \sum_{k=0}^{n-1} \mathbb E_x \exp(-\varphi_{k-1}) L^{c} 
h(X_k). 
\end{equation}
In other words, the process
\begin{equation}\label{mart1}
M_n:= \exp(-\varphi_{n-1})\, h(X_n) - h(x) - \sum_{k=0}^{n-1} \exp(-\varphi_{k-1})L^{c} 
h(X_k), \quad n\ge 0, 
\end{equation}
is a martingale. 
\end{lemma}

~

\noindent
{\em Proof.} Let the initial state $x$ be fixed. 
Let us check the base, \underline{$n=0$}. Note that $\varphi_0 = c(x)$, $\varphi_{-1} = 0$, and $L^c h(x)= \exp(-\phi_0){\cal P}h(x) - h(x)$. So, for $n=0$ the formula (\ref{dyn11}) reads, 
\[
h(x) = 
h(x) + \sum_{k=0}^{-1} \mathbb E_x L^c 
h(X_k), 
\]
which is true due to the  standard convention that $\sum\limits_{k=0}^{-1} \cdots = 0$. 

~

\noindent
Let us check the first step, \underline{$n=1$}: 
\begin{eqnarray*}
\mathbb E_x \exp(-c(x))\, h(X_1) =  h(x) + \sum_{k=0}^{0} \mathbb E_x \exp(-\varphi_{-1})L^{c} 
h(X_k) \equiv h(x) + \exp(-c(x)){\cal P} h(x) - h(x), 
\end{eqnarray*}
or, equivalently, 
\begin{eqnarray*}
\mathbb E_x \exp(-c(x))\, h(X_1) =   \exp(-c(x)){\cal P} h(x),  
\end{eqnarray*}
which is also true. 

~

\noindent
The induction step with a general \underline{$n\ge 1$} follows similarly, using the Markov property. Indeed, assume that the formula (\ref{dyn11}) is true for some $n \ge 0$. Then, for $n+1$ we have, 
\begin{eqnarray*}
& \displaystyle \mathbb E_x \exp(-\varphi_{n})\, h(X_{n+1}) -  h(x) - \sum_{k=0}^{n} \mathbb E_x \exp(-\varphi_{k-1}) L^{c} h(X_k)
 \\\\
& \displaystyle =  \mathbb E_x \exp(-\varphi_{n})\, h(X_{n+1}) 
- \mathbb E_x \exp(-\varphi_{n-1})\, h(X_{n}) 
+ \mathbb E_x \exp(-\varphi_{n-1})\, h(X_{n})
 \\\\
& \displaystyle -  h(x) - \sum_{k=0}^{n-1} \mathbb E_x \exp(-\varphi_{k-1}) L^{c} h(X_k) -  \mathbb E_x \exp(-\varphi_{n-1}) L^{c} h(X_n)
 \\\\
& \displaystyle =  \mathbb E_x \exp(-\varphi_{n})\, h(X_{n+1}) 
- \mathbb E_x \exp(-\varphi_{n-1})\, h(X_{n})  
-  \mathbb E_x \exp(-\varphi_{n-1}) L^{c} h(X_n)
 \\\\
& \displaystyle =   \mathbb E_x \left[\mathbb E_x \left(\exp(-\varphi_{n})\, h(X_{n+1}) 
- \exp(-\varphi_{n-1})\, h(X_{n})  
-   \exp(-\varphi_{n-1}) L^{c} h(X_n) | {\cal F}_n\right)\right]
 \\\\
& \displaystyle =   \mathbb E_x \exp(-\varphi_{n-1})
 \left[\mathbb E_x \left(\exp(-c(X_{n}))\, h(X_{n+1}) 
- h(X_{n})  
- L^{c} h(X_n) | X_n\right)\right] =0, 
\end{eqnarray*}
by definition of $L^c$. This completes the induction step, so the Lemma \ref{dyn11} is proved.

\begin{lemma}[Dynkin's formula 5]\label{dynkin5}
Let $\tau$ be a stopping time with 
$$
\mathbb E_x e^{\alpha\tau} < \infty, \quad \forall \; x\in  S, 
$$
for some $\alpha>0$. 
Then for any function $h$ on $S$ and if $c=\epsilon c_1$ and $\epsilon$ is small enough, 
\[
\mathbb E_x e^{-\varphi_{\tau-1}} h(X_\tau) =  h(x) + \mathbb E_x\sum_{k=0}^{\tau-1}  e^{-\varphi_{k-1}} L^c h (X_k). 
\]
\end{lemma}
Recall that for an irreducible Markov chain with values in a finite state space any hitting time has {\em some} finite exponential moment. This will be used in the next subsection. 

~

\noindent
{\em Proof.} 
We conclude from (\ref{dyn11}), or (\ref{mart1}), due to Doob's theorem about stopped martingales, 
\[
\mathbb E_x e^{-\varphi_{(\tau-1)\wedge n}} h(X_{\tau\wedge n}) =  h(x) + \mathbb E_x\sum_{k=0}^{(\tau-1)\wedge n}  e^{-\varphi_{k-1}} L^c h (X_k). 
\]
Now if $\epsilon$ is small enough, then we may pass to the limit as $n\to\infty$, due to the Lebesgue theorem about a limit under the uniform integrability condition. We have, 
\[
\mathbb E_x e^{-\varphi_{(\tau-1)\wedge n}} h(X_{\tau\wedge n}) \to  \mathbb E_x e^{-\varphi_{\tau-1}} h(X_{\tau}), 
\]
and 
\[
h(x) + \mathbb E_x\sum_{k=0}^{(\tau-1)\wedge n}  e^{-\varphi_{k-1}} L^c h (X_k) 
\to h(x) + \mathbb E_x\sum_{k=0}^{\tau-1}  e^{-\varphi_{k-1}} L^c h (X_k), \quad n\to\infty,
\]
as required. The Lemma \ref{dynkin5} is proved.

\subsection{Poisson equation with a potential with a boundary} 
Recall the equation {\em with the boundary}: 
\begin{equation}\label{pewb}
\exp(-c(x)){\cal P}u(x) - u(x)= -f(x), \; x\in S\setminus \Gamma, \quad u(x) = g(x), \; x\in \Gamma, 
\end{equation}
with a boundary $\Gamma \not = \emptyset$.
The natural candidate for the solution is the function\\
\begin{equation}\label{pe2-solb}
 u(x) := {\mathbb E}_x \left(\sum_{n=0}^{\tau-1}  \exp\left(-\varphi_{n-1}
\right) f(X_n) + 
\exp\left(-\varphi_{\tau-1}
\right)g(X_\tau)\right), 
\end{equation}
$\tau = \inf(n\ge 0:\, X_n\in \Gamma)$. If $x\in \Gamma$, then $\tau = 0$, and we agree that the term $\sum\limits_{k=0}^{-1}$ equals zero. 

\begin{theorem}\label{thm_prov}
If the expectation in (\ref{pe2-solb}) is finite then the function $u(x)$ is a unique solution of the equation (\ref{pewb}). 
\end{theorem}
Recall that $\tau$ does have some exponential moment, so if $c=\epsilon c_1$ as in the 
statement of the Lemma \ref{dynkin5}, and if $\epsilon$ is small enough, then the expression in (\ref{pe2-solb}), {\bf indeed, converges}. 

~

\noindent
{\em The proof} of the Theorem \ref{thm_prov} can be established similarly to the proof of the Theorem \ref{thm_pe1sol}. Firstly, if $x\in \Gamma$, then clearly $\tau = 0$, so that $u(x) = g(x)$.
Secondly, if $x\not \in \Gamma$, then clearly $\tau \ge 1$. Then, due to the Markov property and by splitting the sum,  i.e., taking a sum $\sum\limits_{k=1}^{\tau-1}$ and separately considering the term corresponding to $n=0$ which is just $f(x)$,  we have,  
\begin{eqnarray*}
& \displaystyle u(x) = f(x)  + {\mathbb E}_x \left(\sum_{n=1}^{\tau-1}  \exp\left(-\sum_{k=0}^{n-1} c(X_k)\right) f(X_k) + \exp\left(-\sum_{k=0}^{\tau-1} c(X_k)\right)g(X_\tau)\right)
 \\\\
& \displaystyle = f(x) + {\mathbb E}_x {\mathbb E}_x \left[\sum_{n=1}^{\tau-1}  \exp\left(-\sum_{k=0}^{n-1} c(X_k)\right) f(X_k) + \exp\left(-\sum_{k=0}^{\tau-1} c(X_k)\right) g(X_\tau) | X_1\right]
 \\\\
& \displaystyle = f(x) + {\mathbb E}_x \exp(-c(x)){\mathbb E}_{X_1} \left(\sum_{n=0}^{\tau-1}  \exp\left(-\sum_{k=0}^{n-1} c(X_k)\right) f(X_k) + \exp\left(-\sum_{k=0}^{\tau-1} c(X_k)\right)g(X_\tau)\right)
 \\\\
& \displaystyle = f(x) + \exp(-c(x)) {\mathbb E}_x u(X_1) = f(x) + \exp(-c(x)) {\cal P} u (x), 
\end{eqnarray*}
which shows exactly the equation (\ref{pewb}), as required. The Theorem \ref{thm_prov} is proved.

\subsection{Poisson equation with a potential without a boundary} 
Recall the equation (\ref{eq020}): 
\[
\exp(-c(x)){\cal P}u(x) - u(x)= -f(x), 
\]
and the natural candidate for the solution ``in the whole space'' is the function\\
\begin{equation}\label{pe2-solnb}
 u(x) := \sum_{n=0}^\infty \mathbb E_x \exp\left(-\sum_{k=0}^{n-1} c(X_k)
\right) f(X_k),
\end{equation}
The main question here is the question of convergence. As was mentioned earlier, we are interested in the following case: 
$c(x) = \epsilon c_1(x)$, $\epsilon>0$ small, and $\sum c_1(x) \mu(x) > 0$, where $\mu$ is the unique invariant measure of the Markov chain $X$.

\subsection{Convergence}
The first goal is to justify that \(u\) is well-defined. 
Recall that we want to show convergence of the series, 
\[
 u(x) = \sum_{n=0}^\infty {\mathbb E}_x \exp\left(-\sum_{k=0}^{n-1} c(X_k)\right) f(X_n),
\]
with 
$
c(x) = \epsilon  c_1(x)$, 
$\bar c_1=\sum c_1(x) \mu(x) > 0,$ with $\epsilon >0 $ small.  
Denote 
$$
H_n(\beta, x) := n^{-1} \ln {\mathbb E}_x \exp 
\left(\beta \sum_{k=0}^{n-1} c_1(X_k)\right), \quad \beta \in \mathbb R^1, 
$$ 
or, equivalently, 
\[
 {\mathbb E}_x \exp \left(\beta \sum_{k=0}^{n-1} c_1(X_k)\right) =  {\mathbb E}_x \exp(n\,H_n(\beta,x)). 
\]
(Note that this notation just slightly differs from how the function $H_n$ -- and in the next formula also $H$ -- was defined in the section \ref{sec:ld}: now it is constructed via the ``additive functional'' related to another function $c_1$. Yet, the meaning is similar, so that there is no need to change this standard notation.)
Let 
$$
H(\beta) := \lim_{T\to\infty} H_T(\beta,x), \quad \beta\in 
\mathbb R^1.
$$
As we have seen in the section \ref{sec:ld}, this limit does exist {\bf for all values 
of $\beta$}. (The fact that in the section  \ref{sec:ld} this was shown for another function and under the centering condition for that function is of no importance because the average may be always subtracted.)

~

Also, it may be proved -- {\em left as an exercise to the reader (here some Lemma from \cite{KLS} about estimating the spectral radius may be useful)} -- that  if $\delta>0$ then there exists $n(\delta)$ such that uniformly in $x$
\begin{equation}\label{da}
\sup_{|\beta| \le B} |H(\beta) - H_n(\beta, x)| \le \delta, \quad n\ge n(\delta). 
\end{equation}
Unlike in the section \ref{sec:ld} where it was assumed that $\bar f = 0$, here we compute, 
$$
H'_n(0,x) = n^{-1} {\mathbb E}_x \sum_{k=0}^{n-1} c_1(X_k), 
$$  
where, as usual, the notation $H'_n(0,x)$ stands for $\partial H'_n(\beta,x)/\partial\beta|_{\beta = 0}$. 
Now, due to the Corollary \ref{cor6} it follows, 
$$
\lim_{n\to\infty} n^{-1} {\mathbb   E}_x \sum_{k=0}^{n-1} c_1(X_k) = 
\bar c_1 = \langle
c_1 , \mu \rangle >0.
$$  
This means that in our case 
$$
H'(0) = \bar c_1 > 0,
$$
and that, at least, in some neighbourhood of zero, 
\begin{equation}\label{Hmonotone}
H(\beta)>0, \quad \beta>0, \qquad \& \qquad H(\beta)<0, 
\quad \beta<0.
\end{equation}
Now, convergence of the sum defining $u$ for each $x$ for \(\epsilon>0\) small enough and uniformly in $x$ -- recall that $|S|<\infty$ -- follows from 
(\ref{Hmonotone}). Indeed, choose $\epsilon>0$ so that 
$H(-\epsilon) < 0$ and for a fixed $\delta= - H(-\epsilon)/2$ also choose $n_0$ such that $|H_n(-\epsilon, x) - H(-\epsilon)| < \delta$, for all $n\ge n_0$ and any $x$. We estimate, for $\epsilon$ small and any $x$ (and with $\epsilon$ independent of $x$),  
\begin{eqnarray*}
|u(x)| \le \|f\|_B \, \sum_{n=0}^\infty {\mathbb E}_x \exp\left(-\varepsilon \sum_{k=0}^{n-1} c_1(X_k)\right)
= \|f\|_B \, \sum_{n=0}^\infty\exp(n H_n(-\epsilon,x))
 \\\\
\le  \|f\|_B \, \sum_{n=0}^\infty \exp(n (H(-\epsilon) +  \delta))  
\le  \|f\|_B \, \exp(\delta) \sum_{n=0}^\infty\exp(n H(-\epsilon)) < \infty. 
\end{eqnarray*}

\subsection{$u$ solves the equation}

Let us argue why the function $u$ solves the Poisson equation (\ref{eq020}). By the Markov property, 
\begin{eqnarray*}
u(x) = f(x) +  \exp(-c(x))\sum_{y} {\mathbb E}_x 1(X_1=y) {\mathbb E}_{y} \sum_{k=0}^{\infty} \exp(-\varphi_{k-1}) f(X_k)
 \\\\
= f(x) +  \exp(-c(x))\sum_{y} p_{xy} v(y) =  f(x) + \exp(-c(x)){\mathbb E}_x u(X_1).
\end{eqnarray*}
From this, it follows clearly that, as required, 
$$
 L^cu(x) = \exp(-c(x)){\mathbb E}_x u(X_1) - u(x) = -f(x).
$$

\subsection{Uniqueness of solution}

Uniqueness may be shown in a standard manner. For the difference of two solutions $v = u^1 - 
u^2$ we have $L^c v = 0$. Therefore,  we get, 
\[
v(x) = \exp(-c(x)) {\mathbb E}_x v(X_1). 
\]
After iterating this formula by induction $n$ times, we obtain, 
$$
v(x) = {\mathbb E}_x \exp(- \sum_{k=0}^{n-1}c(X_k)) v(X_n).
$$
Recall that the function $v$ is necessarily bounded on a finite state space $S$. Hence, it follows that $v(x) \equiv 0$. Indeed, we estimate, 
\[
|v(x)| = |{\mathbb E}_x \exp(-\sum_{k=0}^{n-1} c(X_k)) v(X_n)| \le C {\mathbb E}_x \exp(-\sum_{k=0}^{n-1} c(X_k)). 
\]
Hence, we get, for {\em any} $n\ge 0$, 
\begin{eqnarray}\label{eqex}
 |v(x)| \le C {\mathbb E}_x \exp(-\epsilon \sum_{k=0}^{n-1} c_1(X_k)) = C \exp(n H^{}_n(-\epsilon,x)). 
\end{eqnarray}
Recall that $H^{}_n(\beta, x) \to H^{}(\beta), \, n\to\infty$, and that $H^{}(-\epsilon) < 0$ for $\epsilon>0$ small enough. So, the right hand side in (\ref{eqex}) converges to zero exponentially fast with $n\to\infty$. Since the left hand side does not depend on $n$, we get $|v(x)|=0$, i.e., $u^1 \equiv u^2$, as required.

\section{Ergodic Theorem, general case}\label{sec:erg_gen}

Now let us consider a more general construction on a {\bf more general state space}. 
It is assumed that 
\begin{equation}\label{MD}
 \kappa := \inf_{x,x'} \int \left(\frac{P_{x'}(1,dy)}{P_x(1,dy)}\wedge 1 \right)P_x(1,dy) > 0. 
\end{equation}
Note that here $\displaystyle \frac{P_{x'}(1,dy)}{P_x(1,dy)}$ is understood in the sense of the density of the absolute continuous components. For brevity we will be using a simplified notation $P_x(dz)$ for $P_x(1,dz)$. 
Another slightly less general condition will be accepted in the next section but it is convenient to introduce it here: {\bf suppose}  that there exists a measure  $\Lambda$ with respect to which each measure $P_x(1,dz)$ for any $x$ is absolutely continuous,  
\begin{equation}\label{Lambda}
P_x(1,dz) <\!\!\!< \Lambda(dz), \qquad  \forall \, x \in S.
\end{equation}  
Under the assumption (\ref{Lambda})  we have another representation of the constant $\kappa$ from~(\ref{MD}). 
\begin{lemma}\label{newMD}
Under the assumption (\ref{Lambda}), we have the following representation for the constant from (\ref{MD}), 
\begin{equation}\label{betterMD}
 \kappa = \inf_{x,x'} \int \left(\frac{P_{x'}(1,dy)}{\Lambda(dy)}\wedge \frac{P_{x}(1,dy)}{\Lambda(dy)} \right)\Lambda (dy). 
\end{equation}
\end{lemma}

\noindent
{\em Proof. } Firstly, note that clearly the right hand side in (\ref{newMD}) does not depend on any particular measure $\Lambda$, i.e., for any other measure with respect to which both $P_{x'}(1,dy)$ and $P_{x}(1,dy)$ are absolutely continuous the formula (\ref{MD}) gives the same result. Indeed, it follows straightforward from the fact that if, say, $d\Lambda <\!\!\!< d\tilde\Lambda$ and $d\Lambda = f d\tilde\Lambda$, then we get, 
\begin{eqnarray*}
&\displaystyle  \int \left(\frac{P_{x'}(1,dy)}{\Lambda(dy)}\wedge \frac{P_{x}(1,dy)}{\Lambda(dy)} \right)\Lambda (dy) = 
\int \left(\frac{P_{x'}(1,dy)}{f\tilde\Lambda(dy)}\wedge \frac{P_{x}(1,dy)}{f\tilde\Lambda(dy)} \right)f(y) 1(f(y)>0)\tilde\Lambda (dy) 
 \\\\
&\displaystyle  = \int \left(\frac{P_{x'}(1,dy)}{\tilde\Lambda(dy)}\wedge \frac{P_{x}(1,dy)}{\tilde\Lambda(dy)} \right)1(f(y)>0)\tilde\Lambda (dy).
\end{eqnarray*}
However, $P_{x'}(1,dy)<\!\!\!< \Lambda(dy) = f(y) \tilde \Lambda(dy)$, so for any measurable $A$ we have $\int_A P_{x'}(1,dy) 1(f(y)=0) = 0$ and the same for $P_x(1,dy)$, which means that, actually, 
\[
\int \left(\frac{P_{x'}(1,dy)}{\tilde\Lambda(dy)}\wedge \frac{P_{x}(1,dy)}{\tilde\Lambda(dy)} \right)1(f(y)>0)\tilde\Lambda (dy) = \int \left(\frac{P_{x'}(1,dy)}{\tilde\Lambda(dy)}\wedge \frac{P_{x}(1,dy)}{\tilde\Lambda(dy)} \right)\tilde\Lambda (dy).
\]
Respectively, if there are two reference measure $\Lambda$ and, say, $\Lambda'$, then we may take $\tilde\Lambda = \Lambda + \Lambda'$, and the coefficients computed by using each of the two -- $\Lambda$ and $\Lambda'$ -- will be represented via $\tilde\Lambda$ in the same way. 

~

\noindent
Secondly, 
let $\displaystyle f_x(y) = \frac{P_x(1,dy)}{\Lambda (dy)}(y)$. Then, 
\begin{eqnarray*}
&\displaystyle  \kappa = \inf_{x,x'} \int \left(\frac{P_{x'}(1,dy)}{P_x(1,dy)}\wedge \frac{P_{x}(1,dy)}{P_x(1,dy)} \right)P_x(1,dy) 
 \\\\
&\displaystyle  = \inf_{x,x'} \int \left(\frac{P_{x'}(1,dy)}{f_{x}(y)\Lambda(dy)}\wedge \frac{P_{x}(1,dy)}{f_x(y)\Lambda(dy)} \right)f_x(y)\Lambda (dy) 
= \inf_{x,x'} \int \left(\frac{P_{x'}(1,dy)}{\Lambda(dy)}\wedge \frac{P_{x}(1,dy)}{\Lambda(dy)} \right)\Lambda (dy), 
\end{eqnarray*}
as required. The Lemma \ref{newMD} is proved. 

~

Denote 
\[
\kappa(x,x') : = \int \left(\frac{P_{x'}(1,dy)}{P_x(1,dy)}\wedge 1 \right)P_x(1,dy)
\]
Clearly, for any $x,x'\in S$, 
\begin{equation}\label{jkappa}
\kappa(x,x') \ge \kappa. 
\end{equation}

\begin{lemma}\label{jj}
For any $x,x'\in S$,
\[
\kappa(x,x') = \kappa(x',x).
\]
\end{lemma}
\noindent
{\em Proof.} Under the more restrictive assumption (\ref{betterMD}) 
we have, 
\[
\kappa(x',x) = \int \left(\frac{P_{x'}(1,dy)}{P_x(1,dy)}\wedge 1 \right)P_x(dy) 
= \int \left(\frac{P_{x'}(1,dy)}{\Lambda(dy)}\wedge \frac{P_{x}(1,dy)}{\Lambda(dy)} \right)\Lambda(dy), 
\]
which expression is, apparently, symmetric with respect to $x$ and $x'$, as required. Without assuming (\ref{betterMD}) we can argue as follows. Denote $\Lambda_{x,x'}(dz) = P_x(1,dz)+ P_{x'}(1,dz)$. Note that by definition, $\Lambda_{x,x'} = \Lambda_{x',x}$. Then we have, 
\begin{eqnarray}
& \displaystyle \kappa(x',x) = \int \left(\frac{P_{x'}(1,dy)}{P_x(1,dy)}\wedge 1 \right)P_x(1,dy) 
= \int \left(\frac{P_{x'}(1,dy)}{P_x(1,dy)}\wedge 1 \right)\frac{P_x(1,dy)}{\Lambda_{x,x'}(dy)} \Lambda_{x,x'}(dy)
 \nonumber \\ \nonumber \\
& \displaystyle = \int \left(\frac{P_{x'}(1,dy)}{\Lambda_{x,x'}(dy)} \wedge \frac{P_x(1,dy)}{\Lambda_{x,x'}(dy)} \right) \Lambda_{x,x'}(dy). \label{newgamma}
\end{eqnarray}
The latter expression is  symmetric with respect to $x$ and $x'$, which proves the Lemma \ref{jj}. 

\begin{definition}
If an MC $(X_n)$ satisfies the condition (\ref{MD}) -- we call it {\bf MD-condition} in the sequel -- then we call this process {\em Markov--Dobrushin's or {\bf MD-process}}.
\end{definition}

This condition in an easier situation of finite chains was introduced by Markov himself \cite{Markov}; later on, for non-homogeneous Markov processes its analogue was suggested and used by Dobrushin \cite{Dobrushin}. So, we call it Markov--Dobrushin's condition, as already suggested earlier by Seneta. Note that in all cases $\kappa \le 1$. The case $\kappa = 1$ corresponds to the i.i.d. sequence $(X_n)$. In the opposite extreme situation where the transition kernels are singular for different $x$ and $x'$, we have $\kappa = 0$. 
The MD-condition (\ref{MD}) -- as well as (\ref{betterMD}) -- is most useful because it provides an {\bf effective} quantitative upper bound for convergence rate of a Markov chain towards its (unique) invariant measure in total variation metric. 

\begin{theorem}\label{thm_erg2}
Let the assumption (\ref{MD}) 
hold true. Then the process $(X_n)$ is ergodic, i.e., there exists  a limiting probability measure $\mu$, which is stationary and such that (\ref{erg-def}) holds true. Moreover, the uniform bound is satisfied for every $n$, 
\begin{equation}\label{exp_bd3}
 \sup_{x}\sup_{A\in {\cal S}} |P_x(n,A) - \mu(A)| \le (1-\kappa)^{n}.
\end{equation}
\end{theorem}
\noindent
{\em Recall} that the {\bf total variation distance} or metric between two probability measures may be defined as 
\[
\|\mu-\nu\|_{TV} := 2 \sup_A (\mu(A) - \nu(A)). 
\]
Hence, the inequality (\ref{exp_bd3}) may be rewritten 
as 
\begin{equation}\label{exp_bd22}
 \sup_{x} \|P_x(n,\cdot) - \mu(\cdot)\|_{TV} \le 2 (1-\kappa)^{n}, 
\end{equation}

\noindent
{\em Proof.}
\noindent
{\bf 1.} 
Denote for any measurable $A \in \cal S$, 
$$
M^{(n)}(A) := \sup_x P_x(n,A), \quad 
m^{(n)}(A) := \inf_x P_x(n,A).
$$
Due to the Chapman--Kolmogorov equation we have, 
\begin{eqnarray*}
&\displaystyle m^{(n+1)}(A)  = \inf_x  P_x(n+1, A) = \inf_x \int P_x(dz)  P_z(n,A)
 \\\\
&\displaystyle \ge \inf_x \int  P_x(dz) m^{(n)}(A) = m^{(n)}(A). 
\end{eqnarray*}
So, the sequence $(m^{(n)}(A))$ does not decrease. Similarly, $(M^{(n)}(A))$ does not increase. We are going to show the estimate
\begin{equation}\label{mn}
(0\le) \;\; M^{(n)}(A)-m^{(n)}(A) \le (1-\kappa)^{n}.
\end{equation}
In particular, it follows that for any $x, y \in S$ we have, 
\begin{equation}\label{mnxy}
| P_x(n,A)- P_y(n,A)| \le (1-\kappa)^{n}.
\end{equation}
More than that, by virtue of (\ref{mn}) and due to the monotonicity ($M^{(n)}(A)$ decreases, while $m^{(n)}(A)$ increases) both sequences $M^{(n)}(A)$ and $m^{(n)}(A)$ have limits, which limits coincide and are uniform in~$A$:
\begin{equation}\label{mn2}
\lim\limits_{n\to\infty} M^{(n)}(A) = \lim\limits_{n\to\infty} m^{(n)}(A) =: m(A), 
\end{equation}
and 
\begin{equation}\label{mn3}
\sup\limits_A |M^{(n)}(A)-m(A)|\vee \sup\limits_A |m^{(n)}(A)-m(A)| \le (1-\kappa)^{n}.
\end{equation}

~

\noindent
{\bf 2.} Let $x,x' \in  S$, and let $\Lambda_{x,x'}$ be some reference measure for both $P_x(1,dz)$ and $P_{x'}(1,dz)$. 
Again by virtue of Chapman--Kolmogorov's equation we have for any $n>1$ (recall that we accept the notations, $a_+ = a\vee 0 \equiv \max(a,0)$, and $a_- = a\wedge 0 \equiv \min(a,0)$), 

\begin{eqnarray}\label{i3}
&\displaystyle P_x(n,A) - P_{x'}(n,A)   = \int [P_x(1,dz) - P_{x'}(1,dz)] P_z(n-1,A)
 \nonumber  \\ \nonumber \\ \nonumber 
&\displaystyle = \int \left(\frac{P_x(1,dz)}{\Lambda_{x,x'}(dz)} - \frac{P_{x'}(1,dz)}{\Lambda_{x,x'}(dz)}\right) \Lambda_{x,x'}(dz)\,P_z(n-1,A)
 \\ \nonumber \\ \nonumber 
&\displaystyle = \int \left(\frac{P_x(1,dz)}{\Lambda_{x,x'}(dz)} - \frac{P_{x'}(1,dz)}{\Lambda_{x,x'}(dz)}\right)_+\, \Lambda_{x,x'}(dz)\,P_z(n-1,A) 
 \\ \nonumber \\
&\displaystyle + \int \left(\frac{P_x(1,dz)}{\Lambda_{x,x'}(dz)} - \frac{P_{x'}(1,dz)}{\Lambda_{x,x'}(dz)}\right)_-\, \Lambda_{x,x'}(dz)\,P_z(n-1,A).
\end{eqnarray}
Further, we have, 
\begin{eqnarray*}
&\displaystyle \int \left(\frac{P_x(1,dz)}{\Lambda_{x,x'}(dz)} - \frac{P_{x'}(1,dz)}{\Lambda_{x,x'}(dz)}\right)_+\, \Lambda_{x,x'}(dz)\,P_z(n-1,A) 
 \\\\
&\displaystyle \le \int \left(\frac{P_x(1,dz)}{\Lambda_{x,x'}(dz)} - \frac{P_{x'}(1,dz)}{\Lambda_{x,x'}(dz)}\right)_+\, \Lambda_{x,x'}(dz)\,M^{(n-1)}(A), 
\end{eqnarray*}
and similarly, 
\begin{eqnarray*}
&\displaystyle \int \left(\frac{P_x(1,dz)}{\Lambda_{x,x'}(dz)} - \frac{P_{x'}(1,dz)}{\Lambda_{x,x'}(dz)}\right)_- \, \Lambda_{x,x'}(dz)\,P_z(n-1,A) 
 \\\\
&\displaystyle \le \int \left(\frac{P_x(1,dz)}{\Lambda_{x,x'}(dz)} - \frac{P_{x'}(1,dz)}{\Lambda_{x,x'}(dz)}\right)_-\, \Lambda_{x,x'}(dz)\,m^{(n-1)}(A), 
\end{eqnarray*}
On the other hand, 
\begin{eqnarray*}
&\displaystyle \int \left(\frac{P_x(1,dz)}{\Lambda_{x,x'}(dz)} - \frac{P_{x'}(1,dz)}{\Lambda_{x,x'}(dz)}\right)_+ \Lambda_{x,x'}(dz)
+ \int \left(\frac{P_x(1,dz)}{\Lambda_{x,x'}(dz)} - \frac{P_{x'}(1,dz)}{\Lambda_{x,x'}(dz)}\right)_- \Lambda_{x,x'}(dz) 
 \\\\
&\displaystyle = \int \left(\frac{P_x(1,dz)}{\Lambda_{x,x'}(dz)} - \frac{P_{x'}(1,dz)}{\Lambda_{x,x'}(dz)}\right) \Lambda_{x,x'}(dz) = 1-1=0. 
\end{eqnarray*}
Thus, we get,  
\begin{eqnarray*}
&\displaystyle M^{(n)}(A)-m^{(n)}(A) = \sup_x P_x(n,A) - \inf_{x'}P_{x'}(n,A)
 \\\\
&\displaystyle \le \sup_{x,x'} \int \left(\frac{P_x(1,dz)}{\Lambda_{x,x'}(dz)} - \frac{P_{x'}(1,dz)}{\Lambda_{x,x'}(dz)}\right)_+\, \Lambda_{x,x'}(dz)\,(M^{(n-1)}(A)-m^{(n)}(A)). 
\end{eqnarray*}
It remains to notice that (recall that $(a-b)_+ = a - a\wedge b \equiv a - \min(a,b)$)
\begin{eqnarray*}
&\displaystyle \int \left(\frac{P_x(1,dz)}{\Lambda_{x,x'}(dz)} - \frac{P_{x'}(1,dz)}{\Lambda_{x,x'}(dz)}\right)_+\, \Lambda_{x,x'}(dz) 
 \\\\
&\displaystyle = \int \left(\frac{P_x(1,dz)}{\Lambda_{x,x'}(dz)} - \frac{P_{x'}(1,dz)}{\Lambda_{x,x'}(dz)}\wedge \frac{P_{x}(1,dz)}{\Lambda_{x,x'}(dz)}\right)\, \Lambda_{x,x'}(dz) = 1-\kappa(x,x')\le 1-\kappa. 
\end{eqnarray*}
Now the bound (\ref{mn}) follows by induction.

~

\noindent
{\bf 3.} Let us establish the existence of at least one stationary distribution. For any $x\in S$ and any measurable $A$,
\begin{equation}\label{qna}
m^{(n)}(A) \le P_x(n,A) \le M^{(n)}(A).  
\end{equation}
Due to (\ref{mn2}) and (\ref{mn3}), $(P_x(n,A))$ is a Cauchy sequence which converges exponentially fast and uniformly with respect to $A$. 
Denote 
\begin{equation}\label{qa}
q(A):=\lim_{n\to\infty} P_x(n,A). 
\end{equation}
Clearly, due to this uniform convergence, $q(\cdot)\ge 0$, $q(S)=1$, and the function $q$ is additive in $A$. More than that, by virtue of the same uniform convergence in $A$ in (\ref{qa}), the function $q(\cdot)$ is also ``continuous at zero'', i.e. it is, actually, a sigma-additive measure. More than that, 
the uniform convergence implies that 
\begin{equation}\label{qqtv}
\|P_x(n, \cdot) - q(\cdot)\|_{TV} \to 0, \quad n\to\infty. 
\end{equation}

~

\noindent
{\bf 4.} Now, let us show stationarity. 
We have, 
\begin{eqnarray*}
& \displaystyle q(A) = \lim_{n\to\infty} P^{}_{x_0}(n,A) = \lim_{n\to\infty} \int P^{}_{x_0}(n-1,dz)P_z(A)  
 \\\\
& \displaystyle = \int q(dz)P_z(A)  + \lim_{n\to\infty} \int (P^{}_{x_0}(n-1,dz)-q(dz))P_z(A). 
\end{eqnarray*}
Here, in fact, the second term equals zero. Indeed,  
\begin{eqnarray*}
&\displaystyle |\int (P^{}_{x_0}(n-1,dz)-q(dz))P_z(A)| 
\le \int |P^{}_{x_0}(n-1,dz)-q(dz)|P_z(A) 
 \\\\
&\displaystyle 
\le \int |P^{}_{x_0}(n-1,dz)-q(dz)|
= \|P^{}_{x_0}(n-1,\cdot)-q(\cdot)\|_{TV} \to 0, \quad n\to\infty. 
\end{eqnarray*}
Thus, we find that 
\[
q(A) = \int q(dz)P_z(A), 
\]
which is the definition of stationarity. This completes the proof of the Theorem \ref{thm_erg2}.

~

\begin{corollary}\label{cor66}
For any bounded Borel function $f$ and any $0\le s<t$, 
\[
\sup_x \left|{\mathbb E}_x (f(X_t) | X_s) - {\mathbb E}_{inv} f(X_t)\right| 
\equiv \sup_x \left|{\mathbb E}_x (f(X_t) - {\mathbb E}_{inv} f(X_t)| X_s)\right| 
\le C_f (1-\kappa)^{t-s}, 
\]
or, equivalently, 
\[
\sup_x |{\mathbb E}_x (f(X_t) | {\cal F}^X_s) - {\mathbb E}_{inv} f(X_t)| 
\le C_f (1-\kappa)^{t-s}, 
\]
\end{corollary}

\section{Coupling method: general version}\label{sec:coupling_g}
This more general version requires a change of probability space so as to construct coupling. Results themselves in no way pretend to be new: we just suggest a presentation convenient for the author. 
In particular, all newly arising probability spaces on each step (i.e., at each time $n$) are explicitly shown. By ``general'' we do not mean that it is the most general possible: this issue is not addressed here. Just it is more general that in the section \ref{sec:coupling_s}, and it is more involved because of the more complicated probability space, and it provides a better constant in the convergence bound. It turns out that the general version requires a bit of preparation; hence, we start with the section devoted to a couple of random variables, while the case of Markov chains will be considered separately in the next section. 

The following folklore yet important lemma answers the following question: suppose we have two distributions, which are not singular, and the ``common area'' equals some positive constant $\kappa$. Is it possible to realise these two distributions on the same probability space so that the two corresponding random variables {\em coincide} exactly with probability $\kappa$? We call one version of this result ``the lemma about three random variables'', and another one ``the lemma about two random variables''.

\begin{lemma}[``Of three random variables'']\label{otreh}
Let $\xi^{1}$ and $\xi^2$ be two random variables on their (without loss of generality different, and they will be made independent after we take their direct product!) probability spaces $(\Omega^1, {\cal F}^1, \mathbb P^1)$ and $(\Omega^2, {\cal F}^2, \mathbb P^2)$ and with densities $p^1$ and $p^2$ with respect to some reference measure $\Lambda$, correspondingly.  Then, if Markov -- Dobrushin's condition holds true, 
\begin{equation}\label{kappa47}
\kappa := \int \left(p^1(x)\wedge p^2(x)\right) \Lambda(dx) > 0, 
\end{equation}
then there exists one more probability space $(\Omega, {\cal F}, \mathbb P)$ and a random variable on it $\xi^3$ (and $\xi^2$ also lives on $(\Omega, {\cal F}, \mathbb P)$, clearly, with the same distribution) such that 
\begin{equation}\label{eq3rv}
{\cal L}(\xi^3) ={\cal L}(\xi^1), \quad \& \quad \mathbb P(\xi^3 = \xi^2) = \kappa. 
\end{equation}
\end{lemma}
\noindent
Here $\cal L$ denotes the distribution of a random variable under consideration.  
Note that in the case $\kappa = 1$ we have $p^1 = p^2$, so we can just assign $\xi^3:= \xi^2$, and then immediately both assertions of (\ref{eq3rv}) hold. Mention that even if  $\kappa$ were equal to zero (excluded by the assumption (\ref{kappa47})), i.e., the two distributions were singular, we could have posed $\xi^3:= \xi^1$, and again both claims in (\ref{eq3rv}) would have been satisfied trivially. Hence, in the proof below it suffices to assume 
\[
0< \kappa < 1.
\]

~

\noindent
{\em Proof of the Lemma \ref{otreh}.} 
{\bf 1: Construction.}
Let
\[
A_1 := \{x: \; p^1(x) \ge p^2(x)\}, \quad A_2 := \{x: \; p^1(x) < p^2(x)\},
\]
We will need two new independent random variables, 
$\zeta \sim U[0,1]$ (uniformly distributed random variable on $[0,1]$) and $\eta$ with the density 
\[
\displaystyle 
p^\eta(x) := \frac{p^1 - p^1\wedge p^2}{\displaystyle \int (p^1 - p^1\wedge p^2)(y)\Lambda(dy)}(x) 
\equiv \frac{p^1 - p^1\wedge p^2}{\displaystyle \int\limits_{A_1} (p^1 - p^1\wedge p^2)(y)\Lambda(dy)}(x).
\] 
Both $\zeta$ and $\eta$ are assumed to be defined on {\em their own} probability spaces. 
Now let ({\bf on the direct product of all these probability spaces}, i.e., of the probability spaces where the random variables $\xi^1, \xi^2, \zeta, \eta$ are defined)
\[
\xi^3:= \xi^2 1(\frac{p^1}{p^1\vee p^2}(\xi^2)\ge \zeta) + \eta 1(\frac{p^1}{p^1\vee p^2}(\xi^2) < \zeta).
\]
We shall see that $\xi^3$ admits all the desired properties.  Denote 
$$
C := \{\omega: \; \frac{p^1}{p^1\vee p^2}(\xi^2)\ge \zeta\}.
$$
Then $\xi^3$ may be rewritten as 
\begin{equation}\label{xi3}
\xi^3 = \xi^2 1(C) + \eta 1(\bar C).
\end{equation}

~

\noindent
{\bf 2: Verification.} Below $\mathbb P$ is understood as the probability arising on the direct product of the probability spaces mentioned earlier. 
Let
\[
c:= \int_{A_1} (p^1(x)-p^2(x))\Lambda(dx) \equiv 
\int_{A_2} (p^2(x)-p^1(x))\Lambda(dx). 
\]
Due to our assumptions we have, 
\begin{eqnarray*}
& \displaystyle c+ \kappa = \int_{A_1} (p^1(x)-p^2(x))\Lambda(dx) + \int \left(p^1(x)\wedge p^2(x)\right) \Lambda(dx)
 \\\\
& \displaystyle =  \int_{A_1} (p^1(x)-p^2(x))\Lambda(dx) + \int_{A_1} \left(p^1(x)\wedge p^2(x)\right) \Lambda(dx) + \int_{A_2} \left(p^1(x)\wedge p^2(x)\right) \Lambda(dx)
 \\\\
& \displaystyle =  \int_{A_1} p^1(x)\Lambda(dx) + \int_{A_2} p^1(x)  \Lambda(dx) = \int_{A_1 \cup A_2} p^1(x)\Lambda(dx) = 1. 
\end{eqnarray*}
So, 
\[
c = 1-\kappa \in (0,1).
\]
Also, 
\[
p^\eta(x) =  \frac{p^1 - p^1\wedge p^2}{c}(x).
 \] 

~

\noindent
Also notice that 
\[
\mathbb P(C | \xi^2) = \frac{p^1}{p^1\vee p^2}(\xi^2),
\]
and recall that on $C$, $\xi^3 = \xi^2$, while on its complement $\bar C$, 
$\xi^3 = \eta$.
Now, for any bounded Borel measurable function $g$ we have, 
\begin{eqnarray*}
&\displaystyle \mathbb E g(\xi^3) = \mathbb E g(\xi^3) 1(C) + \mathbb E g(\xi^3) 1(\bar C) = \mathbb E g(\xi^2) 1(C) + \mathbb E g(\eta) 1(\bar C)
 \\\\
&\displaystyle =  \mathbb E g(\xi^2) \frac{p^1}{p^1\vee p^2}(\xi^2) 
+ \mathbb E g(\eta)(1-\frac{p^1}{p^1\vee p^2}(\xi^2)) 
 \\\\
&\displaystyle = \mathbb E g(\xi^2) \frac{p^1}{p^1\vee p^2}(\xi^2) 
+ \mathbb E g(\eta) \mathbb E (1-\frac{p^1}{p^1\vee p^2}(\xi^2))
 \\\\
&\displaystyle = \int\limits_{A_1 \cup A_2} g(x) \frac{p^1}{p^1\vee p^2}(x) p^2(x)\Lambda(dx) + 
\int\limits_{(A_1)} g(x) p^\eta(x)\Lambda(dx) \times \int\limits_{(A_2)} (1-\frac{p^1}{p^1\vee p^2}(y))p^2(y)\Lambda(dy)
 \\\\
&\displaystyle = \int\limits_{A_1} g(x)  p^2(x)\Lambda(dx)
+ \int_{A_2} g(x) p^1(x)\Lambda(dx)
+ \int\limits_{A_1} g(x) \frac{p^1 - p^2}{c}(x)\Lambda(dx) \times 
\int\limits_{A_2} (p^2 - p^1)(y)\Lambda(dy)
 \\\\
&\displaystyle = \int\limits_{A_1 \cup A_2} g(x) p^1(x)\Lambda(dx) = \mathbb E g(\xi^1). 
\end{eqnarray*}
Here $(A_1)$ in brackets in $\displaystyle \int\limits_{(A_1)} g(x) p^\eta(x)\Lambda(dx)$ is used with the following meaning: the integral is originally taken over the whole domain, but integration outside the set $A_1$ gives zero; hence, only the integral over this domain remains. 
The established equality $\mathbb E g(\xi^3)= \mathbb E g(\xi^1)$ means that ${\cal L}(\xi^3) = {\cal L}(\xi^1)$, as required. 

~

\noindent
Finally, from the definition of $\xi^3$ it is straightforward that 
\[
\mathbb P(\xi^3 = \xi^2) \ge \mathbb  P(C).
\]
So, 
\begin{eqnarray*}
\mathbb P(\xi^3 = \xi^2) \ge  P(C) = \mathbb E\frac{p^1}{p^1\vee p^2}(\xi^2) = \int \frac{p^1}{p^1\vee p^2}(x)p^2(x)\Lambda(dx) 
 \\\\\
=  \int_{A_1} \frac{p^1}{p^1\vee p^2}(x)p^2(x)\Lambda(dx) 
+ \int_{A_2} \frac{p^1}{p^1\vee p^2}(x)p^2(x)\Lambda(dx)
 \\\\
= \int_{A_1} p^2(x)\Lambda(dx) + \int_{A_2} \frac{p^1}{p^1\vee p^2}(x)p^1(x)\Lambda(dx) = \int (p^1 \wedge p^2) \Lambda(dx) = \kappa. 
\end{eqnarray*}
Let us argue why, actually, 
\[
\mathbb P(\xi^3 = \xi^2) = \mathbb  P(C) = \kappa, 
\]
i.e., why the inequality $\mathbb P(\xi^3 = \xi^2) \ge  \mathbb P(C)$ may not be strict. Indeed, $\mathbb P(\xi^3 = \xi^2) >  \mathbb P(C)$ may only occur if $\mathbb P(\eta 1(\bar C) = \xi^2)>0$ (cf. with (\ref{xi3})), or, equivalently, if $\displaystyle \mathbb P\left(\eta 1(\frac{p^1}{p^1\vee p^2}(\xi^2) < \zeta) = \xi^2\right)>0$. However, 
$$
\omega \in \bar C = \{\omega: \; \frac{p^1}{p^1\vee p^2}(\xi^2) < \zeta\}.
$$
implies $p^1(\xi_2)<p^2(\xi_2)$, that is, $\xi_2 \in A_2$.But on this set the density of $\eta$ equals zero. Hence, $\mathbb P(\xi^3 = \xi^2) > \mathbb  P(C)$ is not possible, which means that, in fact, we have proved that \(\mathbb P(\xi^3 = \xi^2) = \mathbb  P(C) = \kappa\), as required. 
The Lemma \ref{otreh} is proved. 

~

\noindent
Here is another, ``symmetric'' version of the latter lemma. 
\begin{lemma}[``Of two random variables'']\label{odvuh}
Let $\xi^{1}$ and $\xi^2$ be two random variables on their (without loss of generality different, which will be made independent after we take their direct product!) probability spaces $(\Omega^1, {\cal F}^1, \mathbb P^1)$ and $(\Omega^2, {\cal F}^2, \mathbb P^2)$ and with densities $p^1$ and $p^2$ with respect to some reference measure $\Lambda$, correspondingly.  Then, if 
\begin{equation}\label{MDkappa}
\kappa := \int \left(p^1(x)\wedge p^2(x)\right) \Lambda(dx) > 0, 
\end{equation}
then there exists one more probability space $(\Omega, {\cal F}, \mathbb P)$ and two random variables on it $\eta^1, \eta^2$ such that 
\begin{equation}\label{eq2rv}
{\cal L}(\eta^j) ={\cal L}(\xi^j), \; j=1,2, \quad \& \quad \mathbb  P(\eta^1 = \eta^2) = \kappa. 
\end{equation}
\end{lemma}

~

\noindent
{\em Proof of the Lemma \ref{odvuh}.} 
{\bf 1: Construction.} We will need now {\em four} new independent random variables, Bernoulli random variable  
$\gamma$ with $\mathbb P(\gamma=0) = \kappa$ and $\zeta^{0,1,2}$ with the densities 
\begin{eqnarray*}
& \displaystyle p^{\zeta^1}(x) := \frac{p^1 - p^1\wedge p^2}{\displaystyle\int (p^1 - p^1\wedge p^2)(y)\Lambda(dy)}(x), \quad
p^{\zeta^2}(x) := \frac{p^2 - p^1\wedge p^2}{\displaystyle\int (p^2 - p^1\wedge p^2)(y)\Lambda(dy)}(x),
 \\
& \displaystyle p^{\zeta^0}(x) := \frac{ p^1\wedge p^2}{\displaystyle\int (p^1\wedge p^2)(y)\Lambda(dy)}(x).
\end{eqnarray*}
We may assume that they are all defined on their own probability spaces and eventually we consider the {\bf direct product} of these probability spaces denoted as $(\Omega, {\cal F}, \mathbb P)$. As a result, they are all defined on one unique probability space and they are independent there. 
Now, on the same product of all  probability spaces just mentioned, let 
\begin{equation}\label{ety}
\eta^1:= \zeta^0 1(\gamma=0) + \zeta^1 1(\gamma \not=0), \quad \& 
\quad \eta^2:= \zeta^0 1(\gamma=0) + \zeta^2 1(\gamma \not=0).
\end{equation}
We shall see that $\eta^{1,2}$ admit all the desired properties claimed in the Lemma.  

~

\noindent
{\bf 2: Verification.} From (\ref{ety}), clearly, 
\[
\mathbb P(\eta^1=\eta^2) \ge  \mathbb P(\gamma=0) = \kappa. 
\]
Yet, we already saw earlier (in slightly different terms) that this may be only an equality, that is, \(\mathbb P(\eta^1=\eta^2) = \mathbb  P(\gamma=0) = \kappa\).

~

\noindent
Next, since $\gamma$, $\zeta^0$ and $\zeta^{1}$ are independent on $(\Omega, \cal F, \mathbb P)$, for any bounded measurable function $g$ we have, 
\begin{eqnarray*}
&\displaystyle \mathbb E g(\eta^1) 
= \mathbb E g(\eta^1)1(\gamma=0) + \mathbb E g(\eta^1)1(\gamma \not=0) 
 \\\\
&\displaystyle  
= \mathbb E g(\zeta^0)1(\gamma=0) + \mathbb E g(\zeta^1)1(\gamma \not=0) = \mathbb E g(\zeta^0) \mathbb E 1(\gamma=0) + \mathbb E g(\zeta^1) \mathbb E 1(\gamma \not=0) 
 \\\\
&\displaystyle = \kappa \int g(y) p^{\zeta^0}(y)\,\Lambda(dy)  
+(1- \kappa) \int g(y) p^{\zeta^1}(y)\,\Lambda(dy) 
 \\\\
&\displaystyle = \kappa \int \frac{p^1\wedge p^2}{\displaystyle \int (p^1\wedge p^2)\Lambda(dy)}(x) \Lambda(dx) + (1- \kappa) \int g(x) \frac{p^1 - p^1\wedge p^2}{\displaystyle \int (p^1 - p^1\wedge p^2)(y)\Lambda(dy)}(x)\Lambda(dx)
 \\\\
&\displaystyle =  \int p^1\wedge p^2 (x) \Lambda(dx) +  \int g(x) (p^1 - p^1\wedge p^2) (x)\Lambda(dx)
= \int g(y) p^1(y)\,dy = \mathbb E g(\xi^1).  
\end{eqnarray*}
For $\eta^2$ the arguments are similar, so also $\mathbb E g(\eta^2) = \mathbb E g(\xi^2)$. The Lemma \ref{odvuh} is proved. 

\begin{remark}
Note that the extended probability space in the proof of the Lemma \ref{odvuh} has the form, 
\[
(\Omega, {\cal F}, \mathbb P) = (\Omega^1, {\cal F}^1, \mathbb P^1) \times (\Omega^2, {\cal F}^2, \mathbb P^2) \times (\Omega^\gamma, {\cal F}^\gamma, \mathbb P^\gamma) \times \prod_{k=0}^{2} (\Omega^{\zeta_k}, {\cal F}^{\zeta_k}, \mathbb P^{\zeta_k}).
\]
\end{remark}

\section{General coupling method for Markov chains}\label{sec:gcm}
Throughout this section the assumption (\ref{betterMD}) is assumed. 
In this section it is explained how to apply general  coupling method as in the section \ref{sec:coupling_g} to Markov chains in general state spaces $(S, {\cal S})$. Various presentations of this method may be found in \cite{Kalash, Lindvall, Nummelin, Thorisson, Vaserstein},  et al. This section follows the lines from \cite{BV}, which, in turn, is based on \cite{Vaserstein}. Note that in \cite{BV} the state space was $\mathbb R^1$; however, in $\mathbb R^d$ all formulae remain the same. Clearly, this may be further extended to more general state spaces, although, we will not pursue this goal here.

Let us generalize the Lemma \ref{odvuh} to a sequence of random variables and present our coupling
construction for Markov chains based on \cite{Vaserstein}. Assume that the process has a transition density $p(x,y)$ with respect to some reference measure $\Lambda$ and consider two versions $(X^1_n), (X^2_n)$ of the same Markov process with two initial distributions respectively, which also have densities with respect to this $\Lambda$ denoted by $p_{X^1_0}$ and $p_{X^2_0}$ (of course, this does not exclude the case of non-random initial states). Let 
\begin{equation}\label{operator_q}
\kappa(u,v):= \int p(u,t) \wedge p(v,t)\,\Lambda(dt), 
\quad \kappa = \inf\limits_{u,v}\kappa(u,v), 
\end{equation}
and 
\begin{equation}\label{q0}
\kappa(0) :=\int p_{X^1_0}(t) \wedge p_{X^2_0}(t)\,\Lambda(dt).
\end{equation}
It is clear that $0\le \kappa(u,v)\le1$ for all $u,v$. Note that $\kappa(0)$ is not the same as $\kappa_0$ in the previous sections.
We assume that $X^1_0$ and $X^2_0$ have different distributions, so $\kappa(0)<1$. Otherwise we obviously have
$X^1_n\stackrel{d}{=}X^2_n$ (equality in distribution) for all $n$, and the coupling can be made trivially, for example, by letting  $\widetilde X^1_n= \widetilde
X^2_n:=X^1_n$.

Let us introduce a new, vector-valued {\bf Markov process} $\left(\eta^1_n,\eta^2_n,\xi_n,\zeta_n\right)$.  If $\kappa_0=0$ then we set
\begin{equation*}
\eta^1_0:=X^1_0,\; \eta^2_0:=X^2_0,\; \xi_0:=0,\; \zeta_0:=1.
\end{equation*}
Otherwise, if $0<\kappa(0)<1$, then we apply the Lemma \ref{odvuh} to the random variables $X^1_0$ and $X^2_0$ so as to create the random variables 
$\eta^1_0$, $\eta^2_0$, $\xi_0$ and $\zeta_0$ (they correspond to $\eta^1, \eta^2, \xi$, and $\zeta$ in the Lemma). Now, assuming that the random variables $\left(\eta^1_n,\eta^2_n,\xi_n,\zeta_n\right)$ have been determined for some $n$, let us show how to construct them for $n+1$. For this aim, we define the transition probability density $\phi$ with respect to the same measure $\Lambda$ for this (vector-valued) process as follows,
\begin{equation}\label{process_eta}
\phi(x,y):=\phi_1(x,y^1)\phi_2(x,y^2)\phi_3(x,y^3) \phi_4(x,y^4),
\end{equation}
where $x=(x^1,x^2,x^3,x^4)$, $y=(y^1,y^2,y^3,y^4)$, and if
 $0<\kappa(x^1,x^2)<1$, then
\begin{align}
&\displaystyle \phi_1(x,u):=\frac{p(x^1,u)-p(x^1,
u)\wedge p(x^2,u)}{1-\kappa(x^1,x^2)}, \quad
\phi_2(x,u):=\frac{p(x^2,u)-p(x^1,u)\wedge p(x^2,u)}{1-\kappa(x^1,x^2)},
 \nonumber\\\nonumber\\\nonumber
&\displaystyle \phi_3(x,u):=1(x^4=1)\frac{p(x^1,u)\wedge
 p(x^2,u)}{\kappa(x^1,x^2)}+1(x^4=0)p(x^3,u),
 \\\nonumber\\ 
&\displaystyle \phi_4(x,u):=1(x^4=1)\left(\delta_1(u)(1-\kappa(x^1,x^2))+ 
\delta_0(u)\kappa(x^1,x^2)\right) +1(x^4=0)\delta_0(u)\label{phi_4}, 
\end{align}
where $\delta_i(u)$ is the Kronecker symbol, $\delta_i(u) = 1(u=i)$, or, in other words, the delta measure concentrated at state $i$. The case $x^4=0$ signifies coupling already realised at the previous step, and $u=0$ means successful coupling at the transition.  
In the degenerate cases, if $\kappa(x^1,x^2)=0$ (coupling at the transition is impossible), then we may set,  e.g., 
$$
\phi_3(x,u):=1(x^4=1)1(0<u<1) + 1(x^4=0)p(x^3,u),
$$ 
and if $\kappa(x^1,x^2)=1$, then we may set 
$$
\phi_1(x,u)=\phi_2(x,u):=1(0<u<1).
$$ 
In fact, in both degenerate cases $\kappa(x^1,x^2)=0$ or 
$\kappa(x^1,x^2)=1$, the functions $\phi_3(x,u)1(x^4=1)$ (or, respectively, $\phi_1(x,u)$ and $\phi_2(x,u)$) can be defined more or less arbitrarily, only so as to keep the property of conditional independence of the four random variables $\left(\eta^1_{n+1},\eta^2_{n+1},\xi_{n+1},\zeta_{n+1}\right)$ given $\left(\eta^1_n,\eta^2_n,\xi_n,\zeta_n\right)$.

\begin{lemma}\label{lemma:2}
Let the random variables $\widetilde X^1_n$ and $\widetilde X^2_n$, for $n\in\mathbb{Z}_+$ be defined   by the following formulae:
\begin{align*}
\widetilde X^1_n:=\eta^1_n 1(\zeta_n=1)+\xi_n 1(\zeta_n=0), \quad 
\widetilde X^2_n:=\eta^2_n 1(\zeta_n=1)+\xi_n 1(\zeta_n=0).
\end{align*}
Then 
\[
\widetilde X^1_n\stackrel{d}{=}X^1_n, \;\;\widetilde
 X^2_n\stackrel{d}{=}X^2_n, \quad \mbox{for all $n\ge 0$.}
\]
Moreover, 
\[
\widetilde X^1_n=\widetilde X^2_n, \quad \forall \; n\ge
n_0(\omega):=\inf\{k\ge0: \zeta_k=0\}, 
\] 
and
\begin{equation}\label{estimate}
\P(\widetilde X^1_n\neq \widetilde
 X^2_n)\le(1-\kappa(0))\, \E\prod_{i=0}^{n-1}
 (1-\kappa(\eta^1_i,\eta^2_i)).
\end{equation}
Moreover, $\left(\widetilde X^1_n\right)_{n\ge 0}$ and $\left(\widetilde X^2_n\right)_{n\ge 0}$ are both  homogeneous Markov processes, and 
\begin{equation}\label{mpequi}
\left(\widetilde X^1_n\right)_{n\ge 0}\stackrel{d}{=}\left(X^1_n\right)_{n\ge 0}, \quad \& \quad 
\left(\widetilde X^2_n\right)_{n\ge 0}\stackrel{d}{=}\left(X^2_n\right)_{n\ge 0}.
\end{equation}
\end{lemma}

Informally, the processes $\eta^1_n$ and $\eta^2_n$ represent $X^1_n$ and $X^2_n$, correspondingly,
under condition that the coupling was not successful until time $n$, while the process $\xi_n$
represents both $X^1_n$ and $X^2_n$ if the coupling does occur no later than at time $n$. The process $\zeta_n$ represents the moment of coupling: the event  $\zeta_n=0$ is equivalent to the event that coupling occurs no later than at time $n$.
As it follows from \eqref{process_eta} and \eqref{phi_4},
\begin{align*}
&\displaystyle \mathbb P(\zeta_{n+1}=0|\zeta_n=0)=1,
 \\ \\
&\displaystyle \mathbb P(\zeta_{n+1}=0|\zeta_n=1,\eta^1_n=x^1,\eta^2_n=x^2)=\kappa(x^1,x^2).
\end{align*}
Hence, if two processes were coupled at time $n$, then they remain coupled at time $n+1$, and if they were not
coupled, then the coupling occurs with the probability $\kappa(\eta^1_n,\eta^2_n)$.
At each time the probability of coupling at the next step is as large as possible, given the current states.

~

\noindent
{\em For the proof of Lemma \ref{lemma:2}} see \cite{BV}. 

~

From the last lemma  a new version of the exponential bound in the Ergodic Theorem may be derived. 
In general, it {\em may} somehow improve the estimate based on the constant $\kappa$ from Markov--Dobrushin's condition (\ref{MD}) or (\ref{betterMD}). In the remaining paragaphs we do not pursue the most general situation {\bf restricting ourselves again to a simple setting} of $|S|<\infty$. 
Introduce the operator $V$ acting on a (bounded continuous) function $h$ on the space $S \times S$ as follows: for $x=(x^1, x^2)\in  S \times S$ and $X_n := (\tilde X^1_n, \tilde X^2_n)$, 
\begin{equation}\label{V}
Vh(x) := (1-\kappa(x^1,x^2)) \mathbb E_{x^1,x^2}h(X_1) \equiv \exp(\psi(x))\mathbb E_{x^1,x^2}h(X_1),  
\end{equation}
where in the last expression $\psi(x):= \ln (1-\kappa(x^1,x^2))$. The aim is now to find out whether the geometric bound $(1-\kappa)^n$ in (\ref{exp_bd}) under the assumption (\ref{Mar2}) is the optimal one, or it could be further improved, let under some additional assumptions. Let us rewrite the estimate (\ref{estimate}) as follows: 
\begin{equation}\label{l2}
\P(\widetilde X^1_n\neq \widetilde X^2_n)\le 
(1-\kappa(0)) V^n 1(x). 
\end{equation}
Note that by definition (\ref{V}), for the {\bf non-negative} matrix $V$ its sup-norm $\|V\| = \|V\|_{B,B}:=\sup\limits_{|h|_B\le 1} |Vh|_B $ equals $\sup\limits_{x} V1(x)$, where $|h|_B := \max\limits_x |h(x)|$ and $1=1(x)$ is considered as a function on $S\times S$ identically equal to one. Note that $\sup\limits_{x} V1(x)=1-\kappa$. 

\noindent
Now the well-known inequality (see, for example, \cite[\S 8]{KLS}) reads, 
\begin{equation}\label{rkappa}
r(V) \le \|V\| = (1-\kappa). 
\end{equation}
Further, from the Perron--Frobenius Theorem it follows (see, e.g., \cite[(7.4.10)]{FW}), 
\begin{equation}\label{l1}
\lim_n \frac1n \, \ln V^n 1(x)  = \ln r(V). 
\end{equation}
The assertions (\ref{l2}) and (\ref{l1}) together lead to the following result. 
\begin{theorem}\label{lastthm}
Let state space $S$ be finite and let the Markov  condition (\ref{Mar2}) be satisfied. Then
\begin{equation}\label{newrate}
\limsup\limits_{n\to\infty} \frac1n \ln \| P_x(n,\cdot) - \mu(\cdot)\|_{TV} \le \ln r(V).  
\end{equation}
\end{theorem}
\noindent
In other words, for any $\epsilon>0$ and $n$ large enough, 
\begin{equation}\label{newrate2}
\|P_x(n,\cdot) - \mu(\cdot)\|_{TV} \le (r(V)+ \epsilon)^n, 
\end{equation}
which is strictly better than (\ref{exp_bd}) if $r(V)<\|V\| = 1-\kappa$ and $\epsilon>0$ is chosen small enough, i.e., so that $r(V)+ \epsilon < 1-\kappa$. 
It is also likely to be true in more general cases for compact operators $V$ where $r(V)+ \epsilon< 1-\kappa$. However, the full problem  remains open.




\end{document}